\documentclass[preprint,12pt]{elsarticle}

\usepackage{graphicx}
\usepackage{amsfonts}
\usepackage{amssymb}
\usepackage{amsmath}
\usepackage{srcltx} % usar este paquete para el WinEdt
\usepackage{dsfont}
\usepackage{xcolor}
\usepackage{stmaryrd}
\usepackage{hhline}
\usepackage{curve2e}
\usepackage{bxeepic}
\usepackage{amsthm}
\usepackage{epstopdf}

\journal{}

\begin{document}
	
	\begin{frontmatter}
		
		
		
		\title{Decoupling mixed finite elements on hierarchical triangular grids for parabolic problems}%\tnoteref{funding}}
				
		\author[a]{A. Arrar\'as\corref{cor1}}
		\ead{andres.arraras@unavarra.es}
		\author[a]{L. Portero}
		\ead{laura.portero@unavarra.es}
		\cortext[cor1]{Corresponding author.}
		\address[a]{Departamento de Ingenier\'{\i}a Matem\'atica e Inform\'atica, Universidad P\'ublica de Navarra, Edificio de Las Encinas, Campus de Arrosad\'{\i}a, 31006 Pamplona (Spain)}%\\[-2ex]}

		\begin{abstract}
			In this paper, we propose a numerical method for the solution of time-dependent flow problems in mixed form. Such problems can be efficiently approximated on hierarchical grids, obtained from an unstructured coarse triangulation by using a regular refinement process inside each of the initial coarse elements. If these elements are considered as subdomains, we can formulate a non-overlapping domain decomposition method based on the lowest-order Raviart--Thomas elements, properly enhanced with Lagrange multipliers on the boundaries of each subdomain (excluding the Dirichlet edges). A suitable choice of mixed finite element spaces and quadrature rules yields a cell-centered scheme for the pressures with a local 10-point stencil. The resulting system of differential-algebraic equations is integrated in time by the Crank--Nicolson method, which is known to be a stiffly accurate scheme. As a result, we obtain independent subdomain linear systems that can be solved in parallel. The behaviour of the algorithm is illustrated on a variety of numerical experiments.
		\end{abstract}
		
		\begin{keyword}
			cell-centered finite difference \sep domain decomposition \sep hierarchical grid \sep Lagrange multiplier \sep mixed finite element \sep parabolic problem
			\vspace*{0.25cm}
			\MSC[2010] 35K20 \sep 65M06 \sep 65M55 \sep 65M60 \sep 68W10 \sep 76S05
		\end{keyword}
		
	\end{frontmatter}
	
	% \linenumbers
	

	\section{Introduction}\label{sec:introduction}
	
	Let us consider the parabolic initial-boundary value problem
	\begin{subequations}\label{ibvp}
		\begin{align}
		&p_t+\nabla\cdot\mathbf{u}=f&&\hspace*{-2cm}\mbox{in\,\,}\Omega\times(0,t_f],\label{ibvp:a}\\
		&\mathbf{u}=-K\nabla p&&\hspace*{-2cm}\mbox{in\,\,}\Omega\times(0,t_f],\label{ibvp:b}\\
		&p=p_0&&\hspace*{-2cm}\mbox{in\,\,}\Omega\times\{0\},\label{ibvp:c}\\
		&p=g_D&&\hspace*{-2cm}\mbox{on\,\,}\Gamma_D\times(0,t_f],\label{ibvp:d}\\
		&\mathbf{u}\cdot\mathbf{n}=g_N&&\hspace*{-2cm}\mbox{on\,\,}\Gamma_N\times(0,t_f],\label{ibvp:e}
		\end{align}
	\end{subequations}
	where $\Omega\subset\mathbb{R}^2$ is a convex polygonal domain with boundary $\partial\Omega=\overline{\Gamma}_D\cup\overline{\Gamma}_N$ such that $\Gamma_D\cap\Gamma_N=\emptyset$. In this formulation, $K\equiv K(\mathbf{x})\in\mathbb{R}^{2\times 2}$ is a symmetric and positive definite tensor, and $\mathbf{n}$ is the outward unit vector normal to $\partial\Omega$. In the framework of flow in porous media, $p$ represents the fluid pressure, $\mathbf{u}$ denotes the Darcy velocity and $K$ is the rock permeability divided by the fluid viscosity. In this context, it is common for $K$ to be a piecewise continuous tensor, whose discontinuities are typically related to the presence of layers of different material within the flow domain.

The aim of this paper is to design a numerical scheme for approximating the solution of (\ref{ibvp}). The method will be specially conceived to handle complex geometries and discontinuous tensor coefficients while preserving the computational efficiency and ease of implementation. The use of so-called hierarchical grids permits to reach the compromise between both aspects, since they combine the flexibility of unstructured triangulations --to handle irregular boundaries, or interfaces generated by discontinuous coefficients--, with the efficiency provided by regular meshes. In order to construct a hierarchical grid, we consider an unstructured coarse partition of $\Omega$ into a certain number of triangular elements, and subsequently apply a regular refinement process to each of these coarse elements\footnote{The regular triangulation obtained inside each coarse triangle is known as a three-line mesh, since it is composed of parallel lines to the sides of such a triangle.}. The use of this kind of grids has experienced a great development in recent years, due to their potential for solving a considerable class of problems in large scale simulation (see, e.g., \cite{ber:wel:hul:rud:07,gme:rud:14}). Hierarchical grids are also referred to as semi-structured grids, structured multiblock grids or patch-wise structured grids (cf. \cite{rod:gas:lis:12} and references therein).

In the sequel, the coarse elements of our hierarchy of grids will be considered as subdomains of a suitable non-overlapping decomposition of $\Omega$. For the spatial discretization, we consider an expanded mixed finite element (MFE) approximation to (\ref{ibvp}) using the lowest-order Raviart--Thomas spaces. Since these spaces are separately defined for each subdomain, we need to introduce Lagrange multiplier pressures on the subdomain interfaces (and also on the Neumann boundaries) to guarantee the continuity of the normal components of velocity vectors across subdomains. This idea, first proposed in \cite{arb:daw:kee:94,arb:daw:kee:whe:yot:98}, will later permit to decouple the global problem into a set of independent subdomain problems that can be efficiently solved in parallel.

The key idea behind expanded methods is to introduce an auxiliary variable, called the adjusted gradient, as a third explicit unknown (see \cite{arb:daw:kee:94,che:98,gat:heu:01}). In this way, unlike the standard mixed method, the resulting scheme avoids inverting $K$, thus allowing for a non-negative tensor coefficient. Further, it has the additional advantage that numerical quadrature can be used to obtain a cell-centered finite difference scheme for the pressures (cf. \cite{arb:daw:kee:94,arb:daw:kee:whe:yot:98}). Similar ideas leading to cell-centered methods on triangular meshes can be found in \cite{bre:for:mar:06,ewi:sae:she:98,mic:sac:sal:01}. In our case, this strategy yields a stencil-based formulation of the original problem that takes advantage of the structure of the hierarchical grid. In particular, if the material properties (given by $K$) are constant within each subdomain, a unique 10-point stencil will be associated to every interior cell in the subdomain. Moreover, this stencil will be constant regardless of the level of refinement of the hierarchical grid, thus drastically reducing the computational complexity of the problem at hand.

The expanded MFE method, combined with numerical quadrature and enhanced with Lagrange multipliers at the subdomain interfaces, achieves second-order superconvergence for the pressure variable at the centroids of the triangles (see \cite{arb:daw:kee:whe:yot:98}). On equilateral three-line meshes and if $K$ is the identity matrix, second-order superconvergence is also observed for the normal velocity at the midpoints of the edges (cf. \cite{dup:kee:98}). As we will show experimentally, the order of convergence is slightly less than 2 for the normal velocity if a general three-line mesh and/or a general tensor $K$ are considered. The application of a post-processing technique introduced in \cite{dup:kee:98} will further permit to increase the optimal first-order convergence of the velocity vector to superconvergence of order almost 2.

The continuous-in-time semidiscrete problem leads to a system of semi-explicit differential-algebraic equations (DAEs) of index 1. Specifically, the velocity and gradient can be eliminated from the semidiscrete problem via a static condensation procedure, and the resulting scheme is formulated in terms of pressure and Lagrange multiplier unknowns. This formulation comprises a system of ordinary differential equations (ODEs) subject to certain constraints, that is to say, a system of DAEs in which the pressure represents the differential component and the Lagrange multiplier is the algebraic component. In order to guarantee the same order of convergence for both components, a stiffly accurate Runge--Kutta method is used for the time integration (see, e.g., \cite{hai:wan:96}). In particular, we will consider the simplest member in the Lobatto IIIA family, the well-known Crank--Nicolson method, that achieves second-order convergence for both the differential and algebraic variables (cf. \cite{hai:lub:roc:89}). At each time step, the fully discrete scheme consists of a set of linear systems that decouple across subdomains and, therefore, can be efficiently solved in parallel. This parallel solution strategy has two distinct properties: on the one hand, no iterations are required for convergence; on the other, since the use of hierarchical grids implies the same size (in terms of degrees of freedom) for each subdomain problem, a perfectly balanced workload is attained among parallel processors.

The present work extends the ideas introduced in \cite{arb:daw:kee:94,arb:daw:kee:whe:yot:98} for the elliptic problem to the parabolic case. Further, it provides a full derivation of the explicit expressions for the stencil coefficients, which was lacking so far. In doing so, we apply a novel strategy that relies on the three-line structure of each subdomain partition, and permits an efficient implementation of the proposed algorithm. As discussed in \cite{rod:gas:lis:12} (where such an strategy is used for Galerkin finite element discretizations of elliptic problems), the idea is based on a suitable Cartesian distribution of the subdomain nodes.

The solution of evolutionary problems of the form (\ref{ibvp}) has also been tackled in the earlier work \cite{arr:por:14}. In that paper, an overlapping domain decomposition approach for solving (\ref{ibvp}), with $g_N=0$, is designed and analyzed. The method combines an expanded MFE scheme on triangles with a first-order time-splitting technique\footnote{Unlike the method presented here, the scheme in \cite{arr:por:14} makes use of a domain decomposition technique that is not related to a coarse triangulation of $\Omega$. Instead, it is constructed to generate a family of partition-of-unity functions that permits to split the discrete elliptic operator into terms of lower dimension yielding a time-splitting scheme.}. However, its use is restricted to three-line meshes and it is not suitable for dealing with unstructured triangulations and/or piecewise continuous tensor coefficients. The newly proposed method succeeds to circumvent these situations by enhancing the method with Lagrange multipliers; in addition, it provides second-order convergence in time.

To conclude, it is remarkable to note that this enhanced variant of the expanded MFE scheme is closely related to the so-called hybrid MFE method. The hybrid method, first proposed in \cite{arn:bre:85}, is defined by adding Lagrange multipliers on every single edge of the triangulation of $\Omega$ (see also \cite{bof:bre:for:13}). By construction, the present scheme, which restricts the use of Lagrange multipliers to the subdomain interfaces, will be in general much more efficient than the hybrid technique: as long as $K$ is continuous over relatively large subdomains and/or the coarse elements in the hierarchical grid are sufficiently refined, the number of Lagrange multiplier unknowns will be much larger for the hybrid formulation (cf. \cite{arb:daw:kee:whe:yot:98}).

The rest of the paper is organized as follows. In Section \ref{sec:continuous}, we introduce the expanded variational formulation of (\ref{ibvp}). The construction of the expanded MFE spatial discretization on hierarchical grids is discussed in Section \ref{sec:semidiscrete}. This section also includes the matrix formulation of the semidiscrete scheme leading to a cell-centered finite difference method. Section \ref{sec:fully:discrete} is devoted to the time integration of the resulting DAE system. The numerical behaviour of the fully discrete scheme is tested on a variety of numerical experiments in Section \ref{sec:experiments}. The paper ends with some concluding remarks and perspectives. In addition, Appendix A includes a full derivation of the pressure stencil associated to both an upward and a downward triangle in a three-line mesh.

	\section{The expanded variational form}\label{sec:continuous}
	
    For a domain $R\subset\mathbb{R}^2$, let $W^{k,p}(R)$ denote the standard Sobolev space, with $k\in\mathbb{R}$ and $1\leq p\leq\infty$. Let $H^k(R)$ be the Hilbert space $W^{k,2}(R)$. In the sequel, we will mainly use the spaces of square-integrable scalar and vector functions, $L^2(R)$ and $(L^2(R))^2$, respectively. Such spaces are endowed with the inner product $(\cdot,\cdot)_R$ and norm $\|\cdot\|_R$, so that $\|\varphi\|_R=(\varphi,\varphi)^{1/2}_R$. The subscript $R$ will be omitted whenever $R\equiv\Omega$. For a section $S$ of the domain boundary, $\langle\cdot,\cdot\rangle_S$ represents the $L^2(S)$-inner product or duality pairing. We will also use the space
$$
H(\mbox{div};R)=\{\mathbf{v}\in (L^2(R))^2:\nabla\cdot\mathbf{v}\in L^2(R)\}.
$$
    In order to define an expanded formulation, we need to introduce an auxiliary unknown, the so-called adjusted gradient $\tilde{\mathbf{u}}\equiv\tilde{\mathbf{u}}(\mathbf{x},t)$, given by
	$$
	\tilde{\mathbf{u}}=-G^{-1}\nabla p,
	$$
	where $G\equiv G(\mathbf{x})\in\mathbb{R}^{2\times 2}$ is a symmetric and positive definite matrix related to the geometry of $\Omega$. Considering this new variable, the equation (\ref{ibvp:b}) can be rewritten as
	\begin{subequations}
		\begin{align}
		&G\tilde{\mathbf{u}}=-\nabla p&&\hspace*{-3cm}\mbox{in\,\,}\Omega\times(0,T],\label{adjusted:a}\\
		&\mathbf{u}=KG\tilde{\mathbf{u}}&&\hspace*{-3cm}\mbox{in\,\,}\Omega\times(0,T].\label{adjusted:b}
		\end{align}
	\end{subequations}
	The system given by the equations (\ref{ibvp:a}), (\ref{adjusted:a}) and (\ref{adjusted:b}), together with the corresponding initial and boundary conditions, is usually referred to as expanded mixed formulation.
	
	Let us consider that $\Omega$ can be decomposed into a set of $s$ non-overlapping triangular subdomains $\{\Omega_1,\Omega_2,\ldots,\Omega_s\}$, defined to be shape-regular of diameter $H$ and pairwise disjoint, i.e.,
	\begin{equation}\label{decomposition}
	\overline{\Omega}=\cup_{i=1}^s\,\overline{\Omega}_i,\qquad\Omega_{i}\cap\Omega_{j}=\emptyset,\,\,\mathrm{for}\,\,i\neq j.
	\end{equation}
	Let $\partial\Omega_i$ denote the boundary of subdomain $\Omega_i$. This decomposition is assumed to be geometrically conforming in the sense that any edge of $\partial\Omega_i$ is either part of $\partial\Omega$ or coincides with an edge of an adjacent subdomain. If $\Omega_i$ and $\Omega_j$ are adjacent to each other, we denote $\Gamma_{ij}=\partial\Omega_i\cap\partial\Omega_j$, for $i\neq j$, assuming that $\Gamma_{ij}=\Gamma_{ji}$. Furthermore, we set $\Gamma_i=\partial\Omega_i\setminus\Gamma_D$ and define $\Gamma=\cup_{i=1}^s\Gamma_i$. Typically, as we will see below, these subdomains correspond to the elements in a coarse partition $\mathcal{T}_{H}$ of $\Omega$ with mesh size $H$.
	
	Next, we define the subdomain spaces $V_i=H(\mathrm{div};\Omega_i)$, $\tilde{V}_i=(L^2(\Omega_i))^2$ and $W_i=L^2(\Omega_i)$, and the global spaces
	$$
	V=\bigoplus_{i=1}^s\,V_i,\qquad\tilde{V}=\bigoplus_{i=1}^s\,\tilde{V}_i=(L^2(\Omega))^2,\qquad W=\bigoplus_{i=1}^s\,W_i=L^2(\Omega).
	$$
	Note that, in the definition of $V$, we are relaxing the continuity of the space $H(\mathrm{div};\Omega)$ across the subdomain interfaces. Thus, we need to impose it weakly by introducing the additional space $L=H^{1/2}(\Gamma)$. The following weak form can be obtained by integrating the original equations over each $\Omega_i$ and summing up: \emph{Find} $(\mathbf{{u}},\tilde{\mathbf{u}},p,\lambda):[0,T]\rightarrow V\times\tilde{V}\times W\times L$ \emph{such that}
	\begin{subequations}\label{variational}
		\begin{align}
		\hspace*{-0cm}&\displaystyle(p_{t},w)+\sum_{i=1}^s\,(\nabla\cdot\mathbf{u},w)_{\Omega_i}=(f,w)&&\hspace*{0cm}\forall\,{w}\in W,\label{variational:a}\\[-0ex]
		&\displaystyle(G\tilde{\mathbf{u}},\mathbf{{v}})=\sum_{i=1}^s\,\left((p,\nabla\cdot\mathbf{{v}})_{\Omega_i}-\langle\lambda,\mathbf{{v}}\cdot\mathbf{{n}}\rangle_{\Gamma_i}\right)
		-\langle g_D,\mathbf{{v}}\cdot\mathbf{n}\rangle_{\Gamma_D}&&\hspace*{0cm}\forall\,\mathbf{{v}}\in V,\label{variational:b}\\[-0ex]
		&(G\mathbf{{u}},\tilde{\mathbf{v}})=(GKG\tilde{\mathbf{u}},\tilde{\mathbf{v}})&&\hspace*{0cm}\forall\,\tilde{\mathbf{v}}\in\tilde{V},\label{variational:c}\\[-0ex]
		&\displaystyle\sum_{i=1}^s\,\langle\mathbf{u}\cdot\mathbf{{n}},\mu\rangle_{\Gamma_i}=\langle g_N,\mu\rangle_{\Gamma_N}&&\hspace*{0cm}\forall\,\mu\in L,\label{variational:d}\\[-0ex]
		&p(0)=p_0.\label{variational:e}
		\end{align}
	\end{subequations}
	Note that the flux continuity equation (\ref{variational:d}) implies that $\mathbf{u}$ is, in fact, an element of $H(\mathrm{div};\Omega)$. This formulation is referred to as the macro-hybrid expanded variational form of (\ref{ibvp}) with respect to the decomposition (\ref{decomposition}) (cf. \cite{bof:bre:for:13}).

	\section{The expanded mixed finite element method}\label{sec:semidiscrete}

In this section, we describe the spatial discretization of the variational formulation (\ref{variational}) on a hierarchical grid. To this end, we first define the lowest-order Raviart--Thomas spaces for each single subdomain, and then introduce a space of Lagrange multipliers that ensures the flux continuity across subdomains. The subsequent application of a suitable quadrature rule for computing certain vector inner products permits to reduce the expanded MFE method to a cell-centered finite difference scheme for the pressure, enhanced with some additional unknowns on the subdomain interfaces. A static condensation procedure yields the Schur complement form of the semidiscrete system.

	\subsection{A hierarchical triangular grid}
	
\begin{figure}[t]
		\begin{center}\hspace*{-1.8cm}
			\begin{minipage}[t]{0.48\textwidth}
				\begin{center}\includegraphics[scale=0.18]{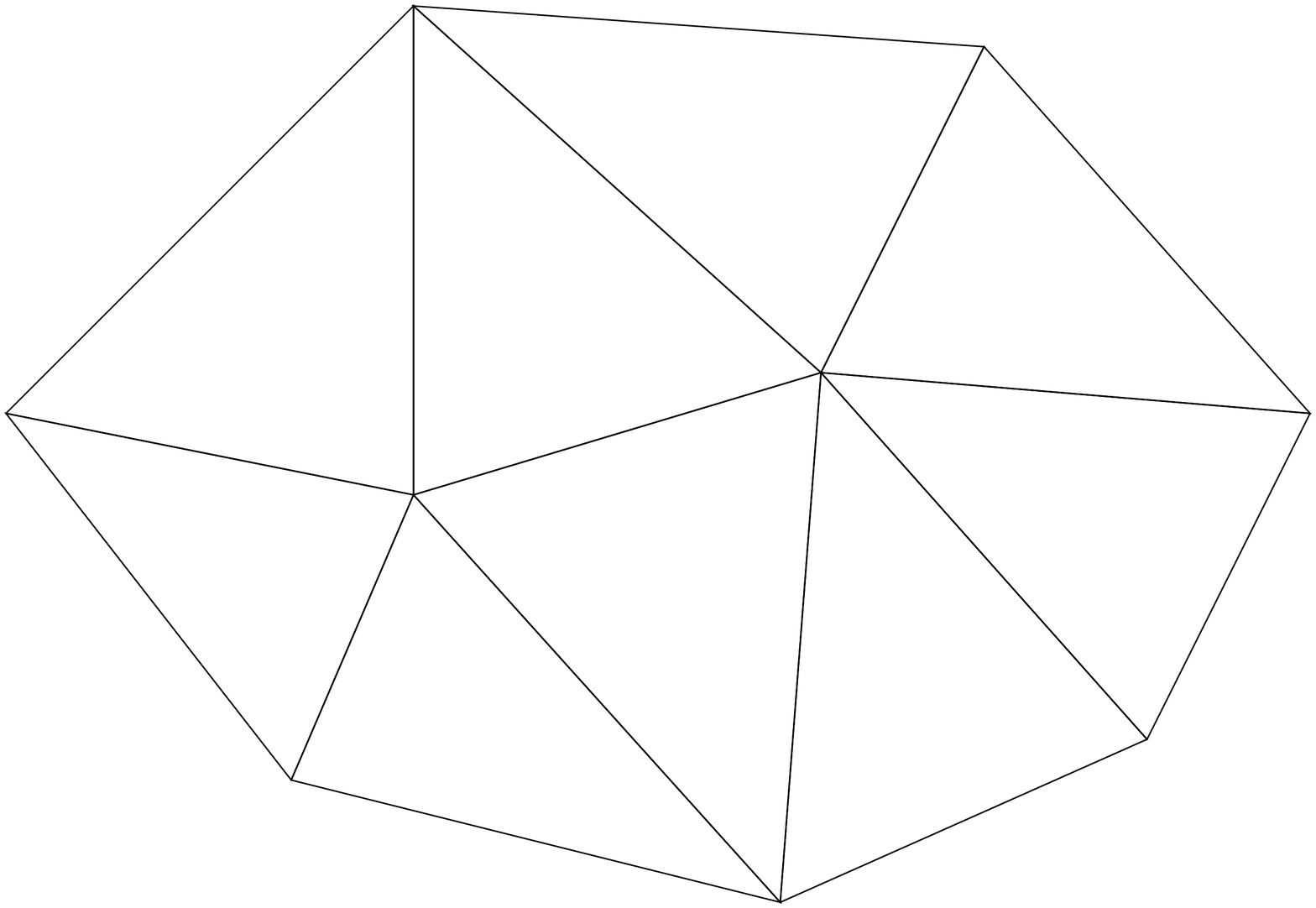}\end{center}
			\end{minipage}
			\hspace*{-0.35cm}
			\begin{minipage}[t]{0.48\textwidth}\begin{center}\includegraphics[scale=0.18]{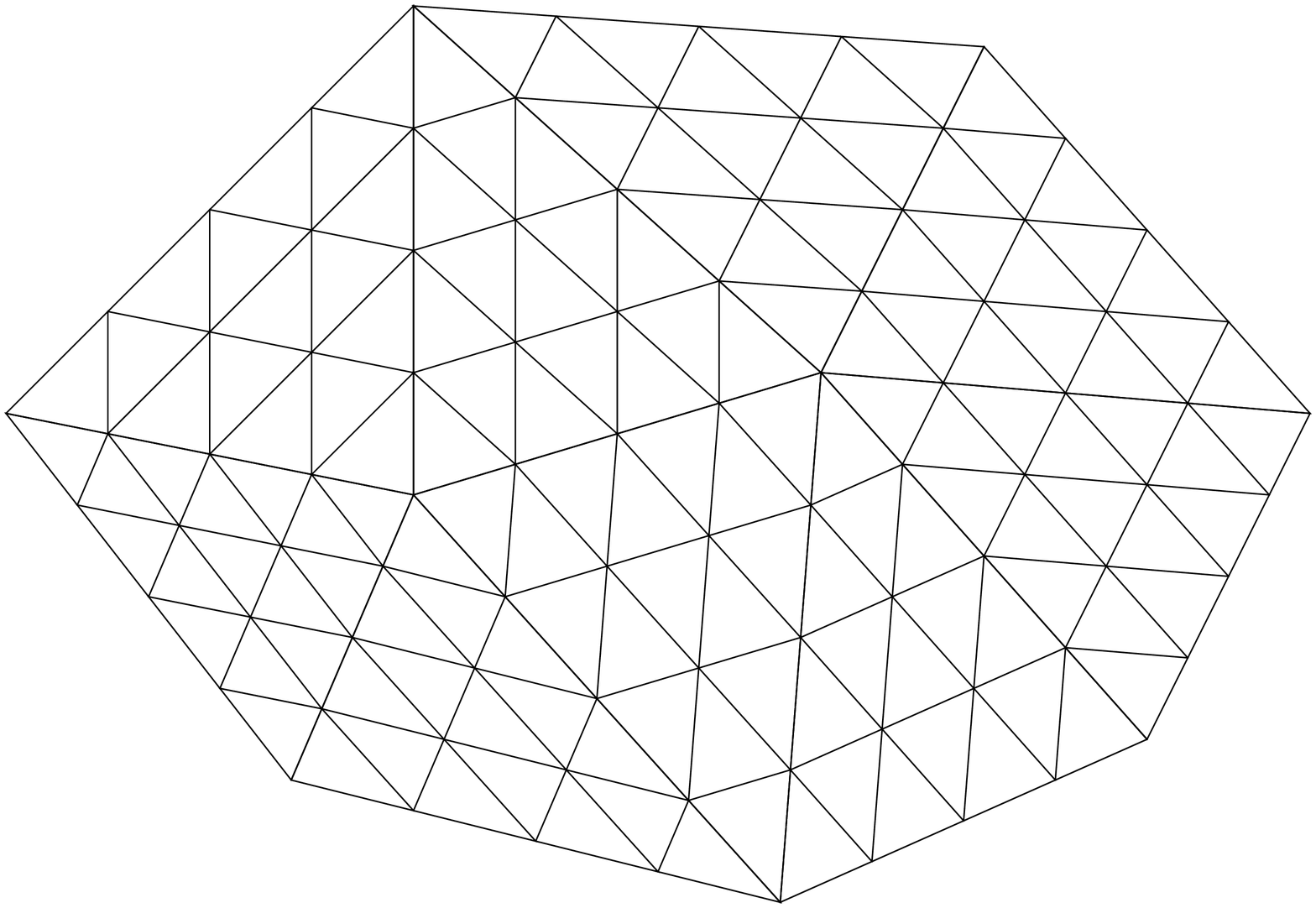}\end{center}
			\end{minipage}
		\end{center}
		\vspace*{-0.5cm}
		\caption{Coarse grid $\mathcal{T}_H$ (left) and fine grid $\mathcal{T}_h\equiv\mathcal{T}_h^2$ (right) for a polygonal domain $\Omega$.}\label{figure:hierarchical:grid}
	\end{figure}

	Let $\mathcal{T}_H$ be an unstructured coarse triangulation of $\Omega$, whose elements are assumed to be the subdomains $\Omega_i$, for $i=1,2,\ldots,s$. This partition is constructed to adequately represent the geometry of the domain, and may also take into account physical features of the problem, such as material properties. If we divide each subdomain $\Omega_i$ into four congruent triangles by connecting the midpoints of its edges, a new regular mesh $\mathcal{T}_{h,i}^1$ is obtained per subdomain. We assume that $\mathcal{T}_{h,i}^1$ and $\mathcal{T}_{h,j}^1$ match on $\Gamma_{ij}$, for $i\neq j$, so that $\mathcal{T}_h^1=\cup_{i=1}^s\mathcal{T}_{h,i}^1$ is a conforming triangulation of $\Omega$. This regular refinement process can be subsequently applied in order to obtain a nested hierarchy of conforming meshes $\mathcal{T}_H\subset\mathcal{T}_h^1\subset\mathcal{T}_h^2\subset\ldots\subset\mathcal{T}_h^k\equiv\mathcal{T}_h$, where $\mathcal{T}_h=\cup_{i=1}^s\mathcal{T}_{h,i}$ is a mesh with the desired fine scale $h$ (cf. \cite{rod:gas:lis:12}). Here, $H=\mathrm{max}_{T\in\mathcal{T}_H}\mathrm{diam}(T)$ and $h=\mathrm{max}_{T\in\mathcal{T}_h}\mathrm{diam}(T)$, where $T$ denotes an element of either $\mathcal{T}_H$ or $\mathcal{T}_h$. Figure \ref{figure:hierarchical:grid} shows an example of the coarse and fine triangulations, $\mathcal{T}_H$ and $\mathcal{T}_h\equiv\mathcal{T}_{h}^2$, for a polygonal domain $\Omega$. Note that, for $i=1,2,\ldots,s$, each $\mathcal{T}_{h,i}$ is a three-line mesh.

	\subsection{Mixed finite element spaces}
	
\begin{figure}[t]
		\begin{center}
			%\hspace*{0.15cm}
			\unitlength=0.035cm
			\begin{picture}(350,110)
			%\allinethickness{0.75pt}
			%Dibujamos el tri\'{a}ngulo de referencia
			\put(50,30){\line(1,0){60}}\put(50,30){\line(50,86){30}}\put(110,30){\line(-50,86){30}}
			%Ponemos nombre a los v\'{e}rtices
			\put(32,20){\scriptsize$(-1,0)$}\put(101,20){\scriptsize$(1,0)$}
			\put(66,85){\scriptsize$(0,\sqrt{3})$}
			%Dibujamos los vectores normales a cada lado
			\put(80,30){\vector(0,-1){15}}\put(95,56){\vector(2,1){14}}\put(65,56){\vector(-2,1){14}}
			%Ponemos nombre a los vectores normales
			\put(85,15){{\footnotesize$\hat{\mathbf{n}}_3$}}
			\put(102,68){\footnotesize{$\hat{\mathbf{n}}_1$}}
			\put(50,68){\footnotesize{$\hat{\mathbf{n}}_2$}}
			%
			%Dibujamos la flecha de la transformaci\'{o}n
			\put(142,53){\vector(1,0){60}}
			\put(166,57){\footnotesize{$F_T$}}
			%\allinethickness{0.75pt}
			%
			%Dibujamos un tri\'{a}ngulo general
			\path(230,45)(310,25)(270,85)(230,45)
			%\path(30,40)(90,20)(104,66)(76,86)(30,40)
			%Ponemos nombre a los v\'{e}rtices
			\put(223.5,37){\footnotesize${\mathbf{r}}_{1}$}\put(311,16){\footnotesize${\mathbf{r}}_{2}$}
			\put(268,90){\footnotesize${\mathbf{r}}_{3}$}
			%Dibujamos los vectores normales a cada lado
			\put(270,35){\vector(-1,-4){4}}\put(290,55){\vector(3,2){13}}
			\put(250,65){\vector(-1,1){11}}
			%Ponemos nombre a los vectores normales
			\put(271,19){{\footnotesize$\mathbf{n}_{3}$}}\put(293,68){\footnotesize{$\mathbf{n}_{1}$}}
			\put(242,77.5){\footnotesize{$\mathbf{n}_{2}$}}
			\end{picture}
		\end{center}\vspace*{-0.3cm}
		\caption{Affine mapping $F_T$ from the reference element $\hat{T}$ onto a generic triangle $T\in\mathcal{T}_h$.}
		\label{fig:mapping:mfe}
	\end{figure}
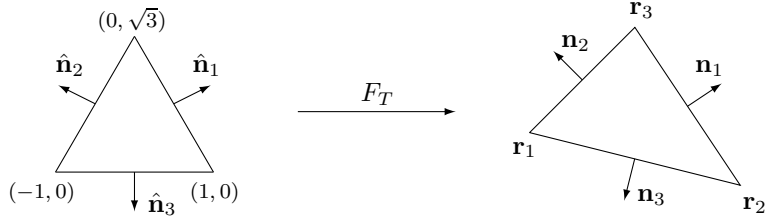

	Let $\hat{T}$ be the reference equilateral triangle with vertices $\hat{\mathbf{r}}_1=(-1,0)^T$, $\hat{\mathbf{r}}_2=(1,0)^T$ and $\hat{\mathbf{r}}_3=(0,\sqrt{3})^T$, and introduce a family of bijective affine mappings $\{F_T\}_{T\in\mathcal{T}_h}$ such that $F_T(\hat{T})=T$. We further define, for each mapping $F_T$, the Jacobian matrix $B_T$ and its determinant $J_T=|\det(B_T)|$. The corresponding vertices of $T$ are denoted by $\mathbf{r}_i=(x_i,y_i)^T$, while the outward unit vectors normal to the edges of $\hat{T}$ and $T$ are represented by $\hat{\mathbf{n}}_i$ and $\mathbf{n}_i$, respectively, for $i=1,2,3$ (see Figure \ref{fig:mapping:mfe}).
	
	Let $\hat{V}_h(\hat{T})$ and $\hat{W}_h(\hat{T})$ be the $RT_0$ finite element spaces on the reference element $\hat{T}$, i.e.,
	$$
	\hat{V}_h(\hat{T})=(\mathbb{P}_{0}(\hat{T}))^2\oplus\hat{\mathbf{x}}\,\mathbb{P}_{0}(\hat{T}),\qquad \hat{W}_h(\hat{T})=\mathbb{P}_{0}(\hat{T}),
	$$
	where $\mathbb{P}_{0}(\hat{T})$ denotes the set of constant functions defined on $\hat{T}$. On every $\hat{e}\subset\partial\hat{T}$, we further define the space $\hat{L}_h(\hat{e})=\mathbb{P}_{0}(\hat{e})$. If $\hat{\mathbf{v}}\in\hat{V}_h(\hat{T})$ and $\hat{w}\in\hat{W}_h(\hat{T})$, the degrees of freedom for $\hat{V}_h(\hat{T})$ are chosen to be the values of $\hat{\mathbf{v}}\cdot\hat{\mathbf{n}}_i$ at the midpoints of the edges, for $i=1,2,3$, while that for $\hat{W}_h(\hat{T})$ is the value of $\hat{w}$ at the element center. Finally, the corresponding degree of freedom for $\hat{\mu}\in\hat{L}_h(\hat{e})$ is the value of $\hat{\mu}$ at the midpoint of the edge $\hat{e}$.
	
	In order to transform any scalar function $\hat{w}$ on $\hat{T}$ or $\hat{\mu}$ on $\hat{e}\subset\partial\hat{T}$ to a generic element $T$ or edge $e\subset\partial T$ belonging to $\mathcal{T}_h$, we introduce the standard isomorphisms
	\begin{align}
	&w\leftrightarrow\hat{w}:&&\hspace*{-3.5cm}w = \hat{w}\circ{F}_T^{-1},\nonumber\\
	&\mu\leftrightarrow\hat{\mu}:&&\hspace*{-3.5cm}\mu = \hat{\mu}\circ{F}_T^{-1}.\nonumber
	\end{align}
	For vector functions $\hat{\mathbf{v}}$ on $\hat{T}$, we use the Piola transformation (cf. \cite{tho:77})
	\begin{align}
	&\mathbf{v}\leftrightarrow\hat{\mathbf{v}}:\quad\mathbf{v}=\left(\frac{1}{J_T}\,{B}_T\,\hat{\mathbf{v}}\right)\circ{F}_T^{-1},\nonumber
	\end{align}
	which is defined to preserve the continuity of the normal components of velocity vectors across interelement edges. This is a necessary condition that must be fulfilled when building approximations to $H(\mathrm{div};\Omega_i)$.

	The subdomain spaces on $\mathcal{T}_{h,i}$, denoted as $V_{h,i}\times\tilde{V}_{h,i}\times W_{h,i}\subset V_{i}\times\tilde{V}_{i}\times W_{i}$, are defined to be
	\begin{align}
	&V_{h,i}=\left\{\mathbf{v}\in V_i: \mathbf{v}|_T\leftrightarrow\hat{\mathbf{v}},\,\hat{\mathbf{v}}\in \hat{V}_h(\hat{T})\,\,\forall\,T\in\mathcal{T}_{h,i}\right\},\nonumber\\[0.5ex]
	&\tilde{V}_{h,i}=\left\{\tilde{\mathbf{v}}\in\tilde{V}_{i}: \tilde{\mathbf{v}}|_T\leftrightarrow\hat{\tilde{\mathbf{v}}},\,\hat{\tilde{\mathbf{v}}}\in \hat{V}_h(\hat{T})\,\,\forall\,T\in\mathcal{T}_{h,i}\right\},\nonumber\\[0.5ex]
	&W_{h,i}=\left\{\,w\in W_{i}: w|_T\leftrightarrow\hat{w},\,\hat{w}\in\hat{W}_h(\hat{T})\,\,\forall\,T\in\mathcal{T}_{h,i}\,\right\}.\nonumber
	\end{align}
	The global spaces on $\mathcal{T}_{h}$ are thus given by
	$$
	V_h=\bigoplus_{i=1}^s\,V_{h,i},\qquad\tilde{V}_h=\bigoplus_{i=1}^s\,\tilde{V}_{h,i},\qquad W_h=\bigoplus_{i=1}^s\,W_{h,i}.
	$$
	Note that the vector functions in $V_h$ have continuous normal components on the edges between elements inside each subdomain, but the space has no continuity constraint on the edges lying on $\Gamma$. In order to recover the continuity for such edges, we introduce the space $L_h\subset L$ of pressure Lagrange multipliers, defined as
	$$
	L_h=\left\{\,\mu\in L:\mu|_e\leftrightarrow\hat{\mu},\,\hat{\mu}\in \hat{L}_h(\hat{e})\,\,\forall\,e\subset(\partial T\cap\Gamma),\,T\in\mathcal{T}_{h}\,\right\},
	$$
	whose elements provide an approximation to $p|_{\Gamma}$.
	
	To determine the dimensions of the preceding spaces, we need to introduce some notations first. Let $N_e$ and $N_T$ be the number of edges and elements in $\mathcal{T}_h$, respectively. Moreover, let $N_e^{\Gamma}$, $N_e^D$ and $N_e^N$ denote the number of edges in $\mathcal{T}_h$ belonging to $\Gamma$, $\Gamma_D$ and $\Gamma_N$, respectively. Then, we have that
	\begin{align*}
	&N_{V}=\mathrm{dim}(V_h)=N_e+N_e^{\Gamma}-N_e^N,\quad
	&&N_{\tilde{V}}=\mathrm{dim}(\tilde{V}_h)=2N_e-N_e^D-N_e^N,\\
	&N_W=\mathrm{dim}(W_h)=N_T,\quad
	&&N_L=\mathrm{dim}(L_h)=N_e^{\Gamma}.
	\end{align*}
	Observe that, in general, $N_e^{\Gamma}\ll N_e$, that is, the number of Lagrange multipliers is much smaller than the total number of edges in $\mathcal{T}_h$. As we will see below, this fact makes the proposed scheme much more efficient than the class of so-called hybrid methods. For such methods, a Lagrange multiplier is introduced on every edge in $\mathcal{T}_h$, so that the dimension of $L_h$ (i.e., the number of Lagrange multiplier unknowns) is $N_e$ instead of $N_e^{\Gamma}$.

	\begin{figure}[t]
		\begin{center}\hspace*{-1.8cm}
			\begin{minipage}[t]{0.48\textwidth}\begin{center}\includegraphics[scale=0.18]{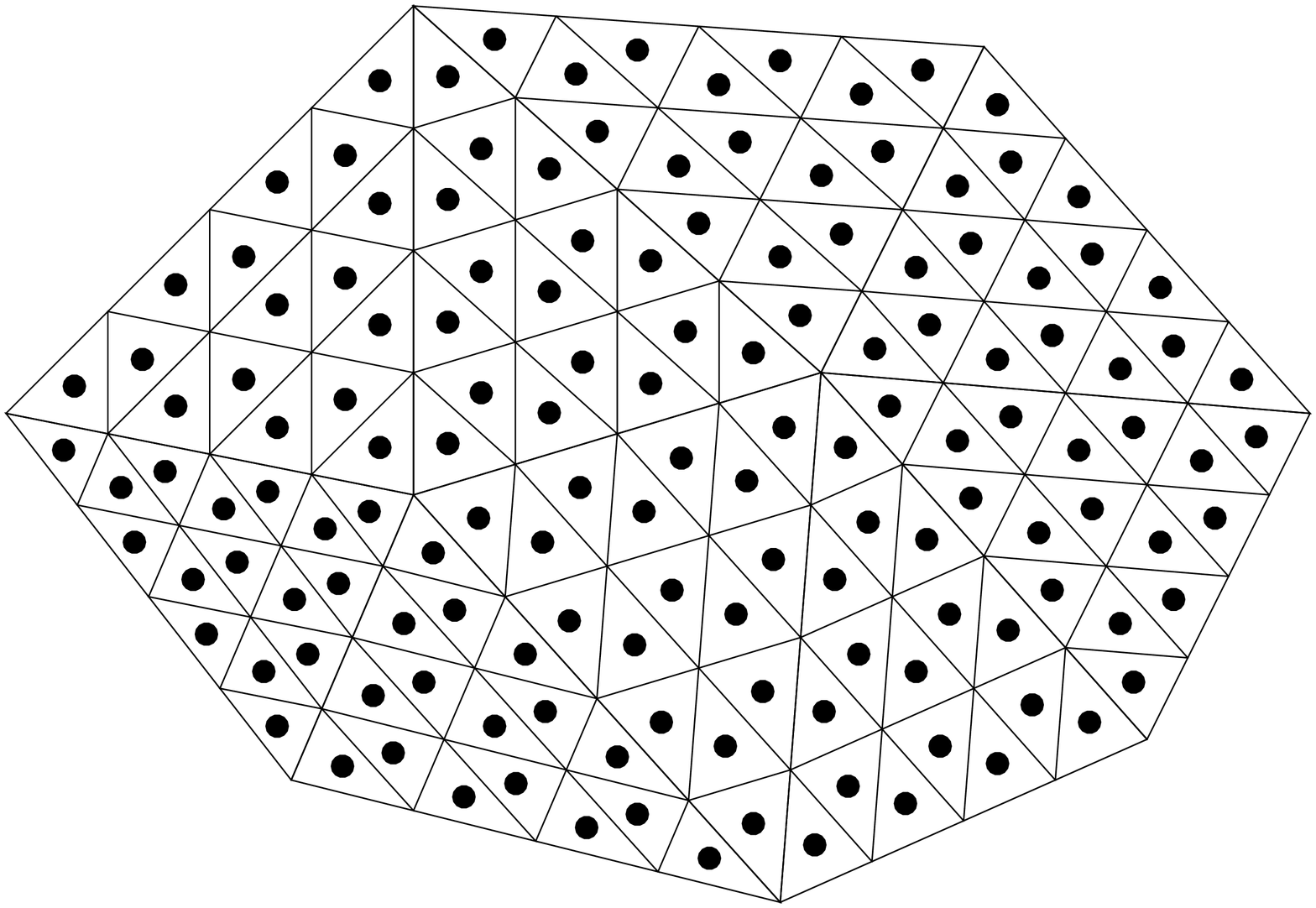}\end{center}
			\end{minipage} \hspace*{-0.35cm}
			\begin{minipage}[t]{0.48\textwidth}\begin{center}\includegraphics[scale=0.18]{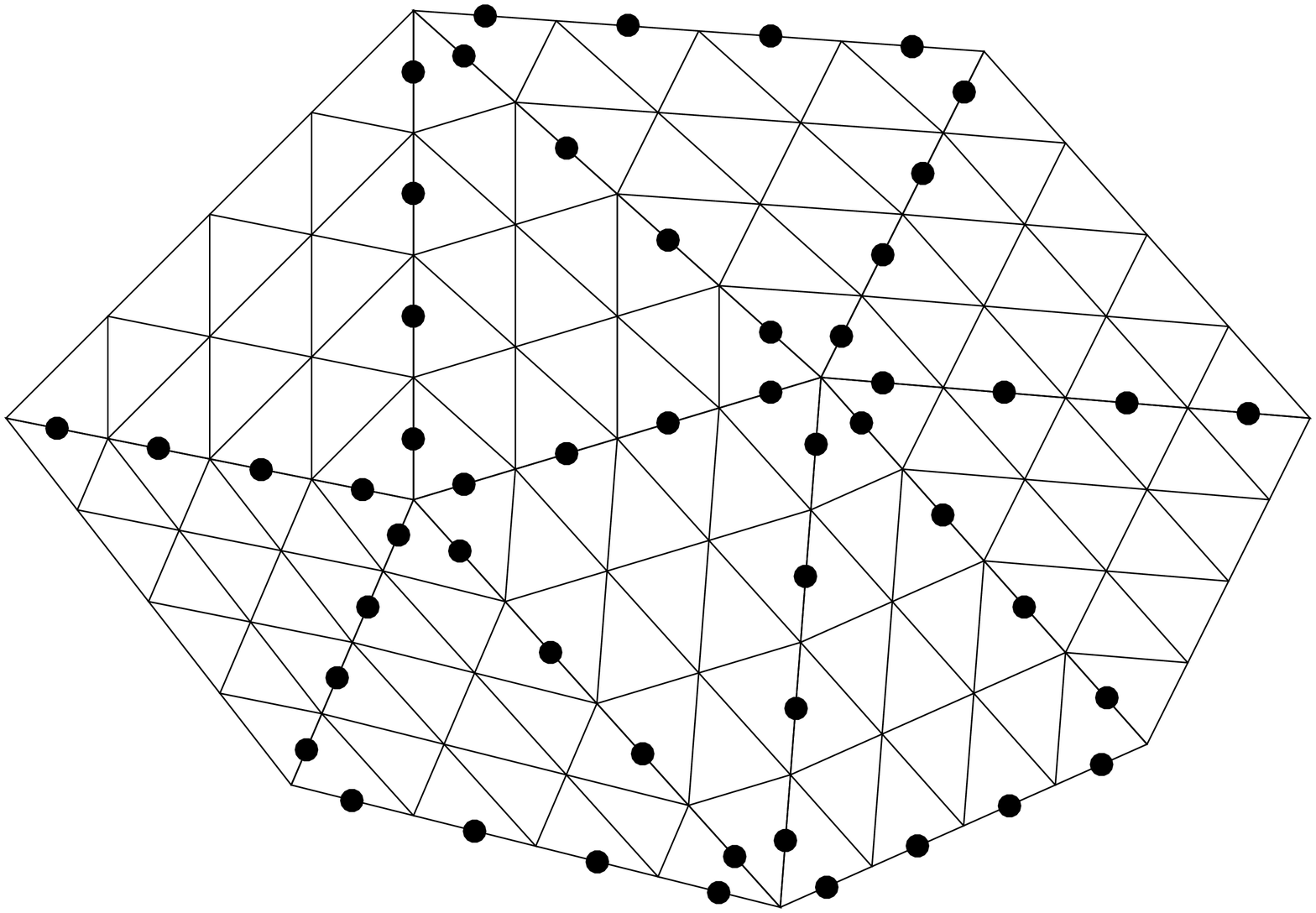}\end{center}
			\end{minipage}
		\end{center}
		\vspace*{-0.5cm}
		\caption{Degrees of freedom for the spaces $W_h$ (left) and $L_h$ (right) on a polygonal domain $\Omega$. The former are located at the centers of the elements $T\in\mathcal{T}_h$, while the latter live on the midpoints of the edges $e\subset(\partial T\cap\Gamma)$. Note that those degrees of freedom in $L_h$ belonging to $\partial\Omega$ correspond to the Neumann boundaries $\Gamma_N$.} \label{figure:dofs}
	\end{figure}

	The degrees of freedom for the pressure spaces $W_h$ and $L_h$ are represented in Figure \ref{figure:dofs} for a polygonal domain $\Omega$. The unknowns $p_h\in W_h$ are located at the centers of the elements $T\in\mathcal{T}_h$, while the Lagrange multipliers $\lambda_h\in L_h$ live on the midpoints of the edges $e\subset(\partial T\cap\Gamma)$. These latter include the edges on the internal boundaries of the subdomains $\Omega_i$, for $i=1,2,\ldots,s$, and the edges on the Neumann boundaries $\Gamma_N$.

	\subsection{The semidiscrete scheme}
	
	The expanded MFE approximation to the variational form (\ref{variational}) reads as follows: \emph{Find} $(\mathbf{{u}}_h,\tilde{\mathbf{u}}_h,p_h,\lambda_h):[0,T]\rightarrow V_h\times\tilde{V}_h\times W_h\times L_h$ \emph{such that}
	\begin{subequations}\label{emfe}
		\begin{align}
		&\hspace*{-0.2cm}\displaystyle(p_{h,t},w)+\sum_{i=1}^s\,(\nabla\cdot\mathbf{u}_h,w)_{\Omega_i}=(f,w)&&\hspace*{-0.15cm}\forall\,{w}\in W_h,\label{emfe:a}\\[0ex]
		&\hspace*{-0.2cm}\displaystyle(G\tilde{\mathbf{u}}_h,\mathbf{{v}})=\sum_{i=1}^s\,((p_h,\nabla\cdot\mathbf{{v}})_{\Omega_i}-\langle \,\lambda_h,\mathbf{{v}}\cdot\mathbf{{n}}\rangle_{\Gamma_i})
		-\langle g_D,\mathbf{{v}}\cdot\mathbf{n}\rangle_{\Gamma_D}&&\hspace*{-0.15cm}\forall\,\mathbf{{v}}\in V_h,\label{emfe:b}\\[0ex]
		&\hspace*{-0.2cm}(G\mathbf{{u}}_h,\tilde{\mathbf{v}})=(GKG\tilde{\mathbf{u}}_h,\tilde{\mathbf{v}})&&\hspace*{-0.15cm}\forall\,\tilde{\mathbf{v}}\in\tilde{V}_h,\label{emfe:c}\\[0ex]
		&\hspace*{-0.2cm}\displaystyle\sum_{i=1}^s\,\langle\mathbf{{u}}_h\cdot\mathbf{{n}},\mu\rangle_{\Gamma_i}=\langle g_N,\mu\rangle_{\Gamma_N}&&\hspace*{-0.15cm}\forall\,\mu\in L_h,\label{emfe:d}\\[0ex]
		&\hspace*{-0.2cm}p_h(0)=\mathcal{S}_hp_0.\label{emfe:e}
		\end{align}
	\end{subequations}
	In this case, the flux matching condition (\ref{emfe:d}) is the discrete analogue to (\ref{variational:d}) and guarantees that $\mathbf{u}_h\in H(\mathrm{div};\Omega)$. The operator $\mathcal{S}_h$ denotes the MFE elliptic projection as defined, e.g., in \cite{tho:97,whe:73}.
	
	Considering the left-hand side of (\ref{emfe:c}) in the previous formulation, we can define a rectangular matrix $S$ whose elements are given by $(S)_{ij}=(G\mathbf{v}_j,\tilde{\mathbf{v}}_i)$, where $\{\mathbf{v}_i\}_{i=1}^{N_V}$ and $\{\tilde{\mathbf{v}}_i\}_{i=1}^{N_{\tilde{V}}}$ are bases for $V_h$ and $\tilde{V}_h$, respectively. In turn, the left-hand side of (\ref{emfe:b}) involves the matrix $S^T$. Following \cite{arb:daw:kee:whe:yot:98}, we subsequently consider $\tilde{V}_h=V_h$ with the aim of making $S$ be a square, symmetric and invertible matrix. In addition, the resulting matrix can be further diagonalized by using a suitable quadrature rule for computing the inner product $(G\mathbf{v}_j,\tilde{\mathbf{v}}_i)$. In this way, the method becomes a cell-centered finite difference scheme for the pressure, with as many unknowns as the number of elements in $\mathcal{T}_h$, $N_T$, enhanced with as many additional unknowns as the number of edges along $\Gamma$, $N_e^{\Gamma}$.
	
	To properly define this quadrature rule, we first obtain the inner product on each $T\in\mathcal{T}_h$ by mapping to the reference element $\hat{T}$. Using the Piola transformation, for any $\mathbf{q}$, $\mathbf{v}\in V_h$ and $\hat{\mathbf{q}}$, $\hat{\mathbf{v}}\in\hat{V}_h(\hat{T})$, we have
	\begin{equation*}
	(G\mathbf{q},\mathbf{v})_T=\left(\frac{1}{J_T}\,B_T^TGB_T\,\hat{\mathbf{q}},
	\hat{\mathbf{v}}\right)_{\hat{T}}.
	\end{equation*}
	If we choose $G|_T=J_TB_T^{-T}B_T^{-1}$ on each element $T$, we simplify the interaction of the functions on $\hat{T}$, thus obtaining $(G\mathbf{q},\mathbf{v})_T=(\hat{\mathbf{q}},\hat{\mathbf{v}})_{\hat{T}}$. Note that $G$ is indeed symmetric and positive definite for every $T\in\mathcal{T}_h$. In virtue of the previous result, the quadrature rule on $T$ is defined as (cf. \cite{arb:daw:kee:94,arb:daw:kee:whe:yot:98})
	\begin{equation*}
	(G\mathbf{q},\mathbf{v})_{Q,T}\equiv
	(\hat{\mathbf{q}},\hat{\mathbf{v}})_{\hat{Q},\hat{T}}=
	\frac{|\hat{T}|}{6}\left(\sum_{i=1}^3\,\hat{\mathbf{q}}(\hat{\mathbf{r}}_{i})\cdot
	\hat{\mathbf{v}}(\hat{\mathbf{r}}_{i})+3\hat{\mathbf{q}}(\hat{\mathbf{x}}_{c})\cdot
	\hat{\mathbf{v}}(\hat{\mathbf{x}}_{c})\right),
	\end{equation*}
	where $|\hat{T}|$ is the area of $\hat{T}$ and $\hat{\mathbf{x}}_c=\frac{1}{3}\,(\hat{\mathbf{r}}_{1}+\hat{\mathbf{r}}_{2}+\hat{\mathbf{r}}_{3})$. This quadrature rule is exact for polynomials of degree $1$. Furthermore, if we denote by $\hat{\mathbf{v}}_i$ the basis function of $\hat{V}_h(\hat{T})$ associated with the $i$-th edge $\hat{e}_i$, for $i=1,2,3$, the following orthogonality condition is satisfied:
	\begin{equation}\label{orthogonality}
	(\hat{\mathbf{v}}_i,\hat{\mathbf{v}}_j)_{\hat{Q},\hat{T}}=
	\begin{cases}
	\,0,&\mathrm{if}\,\,i\neq j,\\
	\,\dfrac{|\hat{T}|}{6}\,|e_i|^2,&\mathrm{if}\,\,i=j,
	\end{cases}
	\end{equation}
	where $|e_i|$ is the length of the corresponding $i$-th edge $e_i$ of $T$, for $i=1,2,3$. The global quadrature rule is thus given by
	\begin{equation}\label{quadrature}
	(G\mathbf{q},\mathbf{v})_{Q}=\sum_{T\in\mathcal{T}_h}(G\mathbf{q},\mathbf{v})_{Q,T}.
	\end{equation}
	Due to (\ref{orthogonality}), the application of this quadrature rule permits to reduce $S$ to a diagonal matrix.
	
	Summarizing, if we choose $\tilde{V}_h=V_h$ and further use the quadrature rule (\ref{quadrature}) for computing the integrals on the left-hand sides of (\ref{emfe:b}) and (\ref{emfe:c}), the semidiscrete scheme is now given by (\ref{emfe:a}), (\ref{emfe:d}) and (\ref{emfe:e}), in combination with the following equations:
	\begin{subequations}\label{semidiscrete}
		\begin{align}
		&\hspace*{-0.25cm}\displaystyle(G\tilde{\mathbf{u}}_h,\mathbf{{v}})_Q=\sum_{i=1}^s\,((p_h,\nabla\cdot\mathbf{{v}})_{\Omega_i}-\langle \,\lambda_h,\mathbf{{v}}\cdot\mathbf{{n}}\rangle_{\Gamma_i})
		-\langle g_D,\mathbf{{v}}\cdot\mathbf{n}\rangle_{\Gamma_D}&&\hspace*{-0.15cm}\forall\,\mathbf{{v}}\in V_h,\label{semidiscrete:b}\\[0ex]
		&\hspace*{-0.25cm}(G\mathbf{{u}}_h,{\mathbf{v}})_Q=(GKG\tilde{\mathbf{u}}_h,{\mathbf{v}})&&\hspace*{-0.15cm}\forall\,{\mathbf{v}}\in V_h.\label{semidiscrete:c}
		\end{align}
	\end{subequations}
	In \cite{arr:por:14}, the unenhanced variant of this method (i.e., that obtained for a single subdomain) was shown to achieve $\mathcal{O}(h)$ optimal convergence for both pressures and velocities. In addition, the pressure variable was experimentally observed to be $\mathcal{O}(h^2)$ superconvergent at the centroids of the triangles. Similar results were derived in \cite{arb:daw:kee:whe:yot:98} for the enhanced mixed method in the elliptic case. As reported in \cite{dup:kee:98}, on equilateral three-line meshes and if $K$ is the identity matrix, the normal velocity also achieves $\mathcal{O}(h^2)$ superconvergence at the midpoints of the edges in the enhanced case. However, the order of convergence for the normal velocity is slightly less than 2 if a general three-line mesh and/or a general tensor $K$ are considered. The application of a post-processing technique further permits to increase the order of convergence of the velocity vector field to almost 2 (cf. \cite{dup:kee:98}). In Section \ref{sec:experiments}, all these results are confirmed through a collection of numerical experiments.

	\subsection{Matrix formulation}\label{subsec:ccfd}
	
	In this subsection, we describe how to express the expanded MFE scheme (\ref{emfe:a}), (\ref{semidiscrete:b}), (\ref{semidiscrete:c}), (\ref{emfe:d}) and (\ref{emfe:e}) in matrix form. Let us denote by $\{\mathbf{v}_i\}_{i=1}^{N_V}$, $\{w_i\}_{i=1}^{N_W}$ and $\{\mu_i\}_{i=1}^{N_L}$ the basis functions of $V_h$, $W_h$ and $L_h$, respectively. Then, the semidiscrete solution $(\mathbf{u}_h,\tilde{\mathbf{u}}_h,p_h,\mathbf{\lambda}_h)$ can be expressed as
	\begin{align*}
	&\mathbf{u}_h(\mathbf{x},t)=\sum_{i=1}^{N_V}U_{h,i}(t)\,\mathbf{v}_i(\mathbf{x}),\quad
	&&\tilde{\mathbf{u}}_h(\mathbf{x},t)=\sum_{i=1}^{N_V}\tilde{U}_{h,i}(t)\,\mathbf{v}_i(\mathbf{x}),\\
	&p_h(\mathbf{x},t)=\sum_{i=1}^{N_W}P_{h,i}(t)\,w_i(\mathbf{x}),\quad
	&&\lambda_h(\mathbf{x},t)=\sum_{i=1}^{N_L}\Lambda_{h,i}(t)\,\mu_i(\mathbf{x}).
	\end{align*}
	Omitting the time dependences, if we further define the vectors
	\begin{align*}
	&U_h=(U_{h,1},U_{h,2},\ldots,U_{h,N_V})^T,\quad
	&&\tilde{U}_h=(\tilde{U}_{h,1},\tilde{U}_{h,2},\ldots,\tilde{U}_{h,N_V})^T,\\
	&P_h=(P_{h,1},P_{h,2},\ldots,P_{h,N_W})^T,\quad
	&&\Lambda_h=(\Lambda_{h,1},\Lambda_{h,2},\ldots,\Lambda_{h,N_L})^T,
	\end{align*}
	the differential system obtained above can be written as
	\begin{equation}\label{full:system}
	\left(
	\begin{array}{c}
	0\\
	0\\
	D{P}'_h \\
	0\\
	\end{array}
	\right)+
	\left(
	\begin{array}{cccc}
	A_1 & -A_2 & 0 & 0\\
	-A_2 & 0 & B^T & -C\\
	0 & B & 0 & 0\\
	0 & -C^{T} & 0 & 0
	\end{array}
	\right)
	\left(
	\begin{array}{c}
	\tilde{U}_h\\
	U_h\\
	P_h\\
	\Lambda_h
	\end{array}
	\right)=
	\left(
	\begin{array}{c}
	0\\
	G^D_h\\
	DF_h \\
	G^N_h\\
	\end{array}
	\right),
	\end{equation}
	where the matrices $A_1,A_2\in\mathbb{R}^{N_V\times N_V}$, $B\in\mathbb{R}^{N_W\times N_V}$ and $C\in\mathbb{R}^{N_V\times N_L}$ are given by
	\begin{align}
	&(A_1)_{ij}=(GKG{\mathbf{v}}_j,{\mathbf{v}}_i),&&\mathrm{for}\,\,i,j\in\{1,2,\ldots,N_V\},\nonumber\\
	&(A_2)_{ij}=(G\mathbf{v}_j,{\mathbf{v}}_i)_Q,&&\mathrm{for}\,\,i,j\in\{1,2,\ldots,N_V\},\nonumber\\
	&(B)_{ij}=(\nabla\cdot\mathbf{v}_j,w_i),&&\mathrm{for}\,\,i=1,2,\ldots,N_W;\,j=1,2,\ldots,N_V,\nonumber\\
	&(C)_{ij}=\langle\mu_j,\mathbf{v}_i\cdot\mathbf{n}\rangle_{\Gamma},&&\mathrm{for}\,\,i=1,2,\ldots,N_V;\,j=1,2,\ldots,N_L.\nonumber
	\end{align}
	By construction, $A_1$ and $A_2$ are block-diagonal matrices with $s$ diagonal blocks, each associated to a subdomain $\Omega_i$, for $i=1,2,\ldots,s$. In particular,
\begin{align*}
&A_1=\mathrm{diag}(A_1^1,A_1^2,\ldots,A_1^s),\\[0ex]
&A_2=\mathrm{diag}(A_2^1,A_2^2,\ldots,A_2^s),
\end{align*}
where $A_1^i,A_2^i\in\mathbb{R}^{(N_V/s)\times (N_V/s)}$. Note that the value $N_V/s\in\mathbb{N}$ corresponds to the number of edges per subdomain. Furthermore, $A_1$ is symmetric, positive definite and sparse, and $A_2$ is diagonal with positive diagonal entries (see (\ref{orthogonality})). In turn, $B$ is a rectangular block matrix of the form
$$
B=\left(\begin{array}{cccc}
	B^1 & & &\\
	& B^2 & &\\
	& & \ddots &\\
	& & & B^{s}
	\end{array}\right),
$$
whose blocks $B^i\in\mathbb{R}^{(N_W/s)\times (N_V/s)}$, for $i=1,2,\ldots,s$. In this case, the value $N_W/s\in\mathbb{N}$ corresponds to the number of elements per subdomain. On the other hand, the diagonal matrix $D\in\mathbb{R}^{N_{W}\times N_{W}}$ is given by
$$
D=\mathrm{diag}(|T_1|,|T_2|,\ldots,|T_{N_W}|),
$$
where $|T_i|$ denotes the area of $T_i$, for $i=1,2,\ldots,N_W$. Finally, the vectors $G_h^D\in\mathbb{R}^{N_V}$, $F_h\in\mathbb{R}^{N_{W}}$ and $G_h^N\in\mathbb{R}^{N_L}$ are defined to be
	\begin{align}
	&(G_h^D)_i=\langle g_D,\mathbf{v}_i\cdot\mathbf{{n}}\rangle_{\Gamma_D},&&\hspace*{-1.5cm}\mathrm{for}\,\,i=1,2,\ldots,N_V,\label{data:a}\\
	&(F_h)_i=\frac{1}{|T_i|}\,(f,w_i),&&\hspace*{-1.5cm}\mathrm{for}\,\,i=1,2,\ldots,N_W,\label{data:b}\\
	&(G_h^N)_i=-\langle g_N,\mu_i\rangle_{\Gamma_N},&&\hspace*{-1.5cm}\mathrm{for}\,\,i=1,2,\ldots,N_L.\nonumber
	\end{align}
	The initial condition $P_h(0)=P_{h}^0\in\mathbb{R}^{N_{W}}$ is obtained as
	\begin{equation}\label{initial:condition}
	(P_{h}^0)_i=\frac{1}{|T_i|}\,\int_{T_i}p_0(\mathbf{x})\,d\mathbf{x},\qquad\mathrm{for}\,\,i=1,2,\ldots,N_W.
	\end{equation}
	The preceding system (\ref{full:system}) can be reduced to a system for ${P}_h$ and $\Lambda_h$ by eliminating both $U_h$ and $\tilde{U}_h$. To this end, we express
	\begin{equation}\label{reduced:system}
	\left(
	\begin{array}{c}
	0\\
	D{P}'_h \\
	0\\
	\end{array}
	\right)+
	\left(
	\begin{array}{cccc}
	\mathcal{A} & \mathcal{B}^T & -\mathcal{C}\\
	\mathcal{B} & 0 & 0\\
	-\mathcal{C}^{T} & 0 & 0
	\end{array}
	\right)
	\left(
	\begin{array}{c}
	\Phi_h\\
	P_h\\
	\Lambda_h
	\end{array}
	\right)=
	\left(
	\begin{array}{c}
	R_h\\
	DF_h \\
	G^N_h\\
	\end{array}
	\right),
	\end{equation}
	where we introduce the submatrix
	\begin{align*}
	&\mathcal{A}=\left(
	\begin{array}{cccc}
	A_1 & -A_2\\
	-A_2 & 0
	\end{array}
	\right),
	\end{align*}
	and the vectors $\mathcal{B}=(0,B)$, $\mathcal{C}=(0,C)^T$, $\Phi_h=(\tilde{U}_h,U_h)^T$ and $R_h=(0,G_h^D)^T$. From the first equation in (\ref{reduced:system}), we obtain
	\begin{align*}
	&\Phi_h=\mathcal{A}^{-1}\left(R_h-\mathcal{B}^TP_h+\mathcal{C}\Lambda_h\right).
	\end{align*}
	Inserting this expression into the second and third equations, the system can be written in the Schur complement form
	\begin{equation}\label{schur:complement}
	\left(
	\begin{array}{c}
	D{P}'_h \\
	0\\
	\end{array}
	\right)+
	\left(
	\begin{array}{cccc}
	M & Q\\
	Q^T & N
	\end{array}
	\right)
	\left(
	\begin{array}{c}
	P_h\\
	\Lambda_h
	\end{array}
	\right)=
	\left(
	\begin{array}{c}
	S_h\\
	T_h\\
	\end{array}
	\right),
	\end{equation}
	where the block elements of the system matrix are given by
	\begin{align*}
	&M=BA_2^{-1}A_1A_2^{-1}B^T,\\
	&N=C^TA_2^{-1}A_1A_2^{-1}C,\\
	&Q=-BA_2^{-1}A_1A_2^{-1}C,
	\end{align*}
	and the right-hand side terms are
	\begin{align*}
	&S_h=DF_h+BA_2^{-1}A_1A_2^{-1}G_h^D,\\
	&T_h=G_h^N-C^TA_2^{-1}A_1A_2^{-1}G_h^D.
	\end{align*}
	Note that $M\in\mathbb{R}^{N_W\times N_W}$ and $N\in\mathbb{R}^{N_V\times N_V}$ are both symmetric and positive definite. Moreover, $M$ is a block-diagonal matrix of the form
$$
M=\mathrm{diag}(M^1,M^2,\ldots,M^s),
$$
where each block $M^i\in\mathbb{R}^{(N_W/s)\times(N_W/s)}$ is associated to a subdomain $\Omega_i$, for $i=1,2,\ldots,s$, and can be obtained as
$$
M^i=B^i(A_2^i)^{-1}A_1^i(A_2^i)^{-1}(B^i)^T.
$$
As a consequence, $M$ decouples across subdomains. Moreover, each block $M^i$ is a sparse matrix with, at most, 10 non-zero entries per row. The specific expressions for the coefficients of the 10-point stencil corresponding to $MP_h$ are derived in detail in Appendix A. Figures \ref{fig:stencil:up} and \ref{fig:stencil:down} further show the stencil of a pressure unknown associated to an upward and a downward triangle, respectively.

	\section{The fully discrete scheme}\label{sec:fully:discrete}
	
	The resulting system (\ref{schur:complement}) can be formulated in the form
	\begin{subequations}\label{semiexplicit:dae}
		\begin{align}
		P'_h(t)&=F_1\left(t,P_h,\Lambda_h\right),\label{semiexplicit:dae:a}\\
		0&=F_2\left(t,P_h,\Lambda_h\right),\label{semiexplicit:dae:b}
		\end{align}
	\end{subequations}
	where
	\begin{align*}
	&F_1\left(t,P_h,\Lambda_h\right)=D^{-1}\left(S_h(t)-MP_h-Q\Lambda_h\right),\\
	&F_2\left(t,P_h,\Lambda_h\right)=T_h(t)-Q^TP_h-N\Lambda_h.
	\end{align*}
	This is a system of semi-explicit DAEs, namely: the system of ODEs (\ref{semiexplicit:dae:a}) subject to the constraint (\ref{semiexplicit:dae:b}). In this context, $P_h(t)$ is usually referred to as the differential component, whereas $\Lambda_h(t)$ is called the algebraic component. Since $\det(N)\neq 0$, (\ref{semiexplicit:dae}) is indeed an index-1 DAE system. The initial condition for $P_h(t)$ is provided by (\ref{initial:condition}), and that for $\Lambda_h(t)$ can be derived from (\ref{semiexplicit:dae:b}) as
	$$
	\Lambda_h(0)=\Lambda_h^0=N^{-1}(T_h(0)-Q^TP_h^0).
	$$
	The initial values are thus consistent with (\ref{semiexplicit:dae}), i.e., $F_2\left(0,P_h^0,\Lambda_h^0\right)=0$.
	
	It is well known that the time integration of index-1 DAE systems requires the application of stiffly accurate Runge--Kutta methods\footnote{A Runge--Kutta method is called stiffly accurate if its coefficients in the Butcher tableau satisfy $a_{qj}=b_j$, for $j=1,2,\ldots,q$.} in order to obtain the same order of convergence for both the differential and algebraic components (cf. \cite{hai:wan:96}). For instance, the families of Lobatto IIIA and Lobatto IIIC methods are convergent of order $2q-2$ for both components, $q$ being the number of internal stages of the corresponding method (cf. \cite{hai:lub:roc:89}). The simplest member in the Lobatto IIIA family (i.e., that corresponding to $q=2$) is the Crank--Nicolson method. Applied to (\ref{semiexplicit:dae}), this method reads
	\begin{equation}\label{fully:discrete}
	\begin{cases}
	\,P_h^0=P_h(0),\\[0.5ex]
	\,\Lambda_h^0=\Lambda_h(0),\\[1ex]
	\,\dfrac{P_h^{n+1}-P_h^{n}}{\tau}+\dfrac{1}{2}\,D^{-1}M\left(P_h^n+P_h^{n+1}\right)+\dfrac{1}{2}\,D^{-1}Q\left(\Lambda_h^n+\Lambda_h^{n+1}\right)\\[1.5ex]
	\hspace*{8cm}=\dfrac{1}{2}\,D^{-1}\left(S_h^n+S_h^{n+1}\right),\\[1.5ex]
	\,Q^TP_h^{n+1}+N\Lambda_h^{n+1}=T_h^{n+1},\qquad n=0,1,\ldots,n_f-1,
	\end{cases}
	\end{equation}
	yielding approximations $P_h^n\approx P_h(t_n)$ and $\Lambda_h^n\approx\Lambda_h(t_n)$ on the equidistant time grid $0=t_0<t_1<\ldots<t_{n_f}=t_f$, where $t_n=n\tau$ and $\tau=t_f/n_f$, for $n=0,1,\ldots,n_f\in\mathbb{N}$. The notations $S_h^{n}=S_h(t_{n})$ and $T_h^n=T_h(t_n)$ are also used.

In order to solve the system (\ref{fully:discrete}), we define $\hat{M}=I+\frac{1}{2}\,\tau D^{-1}M$ and consider the following procedure:
\begin{enumerate}
\item Set $P_h^0$ and $\Lambda_h^0$.
\item For $n=0,1,\ldots,n_f-1$:\\[-2ex]
\begin{enumerate}
\item Solve the system for $\Lambda_h^{n+1}$:
	\begin{align}\label{eq:lambda}
	&\hspace*{-0.8cm}\textstyle\left(N-\frac{1}{2}\,\tau Q^T\hat{M}^{-1}D^{-1}Q\right)\Lambda_h^{n+1}=T_h^{n+1}\nonumber\\
 &\textstyle+Q^T\hat{M}^{-1}\left(\frac{1}{2}\,\tau D^{-1}\left(S_h^n+S_h^{n+1}+Q\Lambda_h^n\right) -\left(I-\frac{1}{2}\,\tau D^{-1}M\right)P_h^n\right).\nonumber\\[1ex]
	\end{align}
\item Solve the system for $P_h^{n+1}$:
	\begin{align}\label{eq:p}
	&\hspace*{-0.2cm}\textstyle\hat{M}P_h^{n+1}=\frac{1}{2}\,\tau D^{-1}\left(S_h^n+S_h^{n+1}\right)-\frac{1}{2}\,\tau D^{-1}Q\left(\Lambda_h^n+\Lambda_h^{n+1}\right)\nonumber\\
&\textstyle\hspace*{6.5cm}+\left(I-\frac{1}{2}\,\tau D^{-1}M\right)P_h^n.
	\end{align}
\end{enumerate}
\end{enumerate}
The solution of (\ref{eq:lambda}) may be obtained by using iterative linear solvers that require the computation of matrix-vector products involving the system matrix $N-\frac{1}{2}\,\tau Q^T\hat{M}^{-1}D^{-1}Q$. In such a case, we need to solve problems of the form
	$$
	\textstyle\hat{M}x=b,
	$$
	for certain right-hand sides $b$. Taking into account that the coefficient matrix $\hat{M}$ decouples across subdomains, this equation involves independent subdomain problems which can be solved simultaneously. In addition, the use of hierarchical grids entails the same size for such subdomain problems, thus yielding a perfectly balanced workload among parallel processors. Once $\Lambda_h^{n+1}$ has been computed from (\ref{eq:lambda}), the same ideas can be applied to the solution of (\ref{eq:p}).

\section{Numerical experiments}\label{sec:experiments}

In this section, we study the numerical behaviour of the proposed method in the solution of a collection of parabolic initial-boundary value problems of the form (\ref{ibvp}): on the one hand, we examine its convergence properties by considering problems with a known analytical solution; on the other, we analyze its qualitative performance when applied to non-stationary flow models in porous media.

\subsection{Convergence examples with known analytical solutions}

\subsubsection{A smooth solution test}

Let us consider (\ref{ibvp}) posed on the irregular polygon with 7 sides shown in Figure \ref{figure:hierarchical:grid}. The polygon vertices are located at the points $(0,0)^T$, $(1,1)^T$, $(2.4,0.9)^T$, $(3.2,0)^T$, $(2.8,-0.8)^T$, $(1.9,-1.2)^T$ and $(0.7,-0.9)^T$. We further consider $t_f=2$, $\Gamma_D=\partial\Omega$, and $$K=\left(\begin{array}{cc}2&1\\1&2\end{array}\right).$$
The functions $f$, $p_0$ and $g_D$ are defined in such a way that
$$p(x,y,t)=\sin(\pi t)\sin(\pi x)\sin(\pi y)$$
is the exact solution of the problem. The spatial domain $\Omega$ is decomposed into $s=9$ subdomains, as shown in Figure \ref{figure:hierarchical:grid} (left).

In the sequel, the pressure and velocity errors are computed by combining the $\ell^{\infty}$-norm in time with various norms in space. In particular, with an abuse of notation, the pressure errors are obtained in the norms
\begin{equation}\label{norm:press:2}
\left\|r_hp-P_h\right\|_{\ell^{\infty}(\ell^2)}=\max_{1\leq n\leq n_f}\left\|r_hp(t_{n})-P_h^{n}\right\|_{\ell^2},
\end{equation}
and $\|r_hp-P_h\|_{\ell^{\infty}(\ell^{\infty})}$, where the maximum norm is used in both time and space. In these expressions, for any given $t$, $r_hp(t)\in\mathbb{R}^{N_W}$ is a vector whose $i$-th component is $p(\textbf{x}_c^i,t)$, $\textbf{x}_c^i$ being the centroid of the $i$-th triangle of $\mathcal{T}_h$, for $i=1,2,\ldots,N_W$. Regarding the velocity variable, we consider the following norms
\begin{align}
&\|\mathbf{u}-\mathbf{R}(\mathbf{u}_h)\|_{\ell^{\infty}(L^2)}= \max_{1\leq n\leq n_f}\|\mathbf{u}(t_{n})-\mathbf{R}(\mathbf{u}_h^n)\|,\label{norm:postproc:flux}\\
&\|\left(\mathbf{u}-\mathbf{u}_h\right)\cdot\mathbf{n}\|_{\ell^{\infty}(\ell^2)}=\max_{1\leq n\leq n_f}\| \left(\mathbf{u}(t_n)-\mathbf{u}_h^n\right)\cdot\mathbf{n}\|_{\ell^2}.\label{norm:normal:flux}
\end{align}
In both expressions, $\mathbf{u}_h^n$ is a function of $V_h$ defined as
$
\mathbf{u}_h^n(\mathbf{x})=\sum_{i=1}^{N_V}U_{h,i}^n\,\mathbf{v}_i(\mathbf{x}),
$
where $U_{h,i}^n$ are the elements of a vector $U_h^{n}$ containing the normal components of the numerical flux, given by
$$
U_h^{n}=A_2^{-1}A_1A_2^{-1}\left( G_h^{D}(t_{n})+B^TP_h^{n}-C\Lambda_h^{n}\right).
$$
In addition, the term $\left(\mathbf{u}(t_n)-\mathbf{u}_h^n\right)\cdot\mathbf{n}$ on the right-hand side of (\ref{norm:normal:flux}) denotes a vector in $\mathbb{R}^{N_V}$, whose elements are the normal components of $\mathbf{u}(t_n)-\mathbf{u}_h^n$ at the midpoints of the edges of $\mathcal{T}_h$. The expression (\ref{norm:postproc:flux}) further includes a linear operator $\mathbf{R}:V_h\rightarrow L^2(\Omega)$ that provides, for any $\mathbf{u}_h^{n}\in V_h$, a post-processed flux $\mathbf{R}(\mathbf{u}_h^{n})$ as defined in \cite{dup:kee:98}. Specifically, $\mathbf{R}(\mathbf{u}_h^{n})$ is a piecewise function given by
$$
\mathbf{R}(\mathbf{u}_h^n)=\left(\begin{array}{c}a\,x+b\,y+c\\d\,x+e\,y+f\end{array}\right)
$$
on each interior triangle $T\in\mathcal{T}_h$. The coefficients $a,b,\ldots,f$ are obtained as the solution of an overdetermined system with 9 equations by a least squares procedure. Such equations are obtained by imposing that the normal components of both the numerical and post-processed fluxes coincide at the midpoints of 9 edges: the three edges of $T$ and the two other edges of each of the three triangles bordering $T$. As described in \cite{dup:kee:98}, for those triangles containing an edge on the boundary of $\Omega$, a suitable modification of this procedure can be applied in order to preserve the global accuracy. In that work, the authors show theoretically that, for MFE methods, this post-processing technique is second-order accurate on three-line meshes if $K$ is the identity matrix. For enhanced MFE methods, numerical evidence reveals that the post-processed flux achieves close to second-order accuracy, even if a full tensor $K$ is considered.

The integral on the right-hand side of (\ref{norm:postproc:flux}) is approximated element-wise by the midpoint quadrature rule. In turn, we use the formula
\begin{equation*}
\int_{T}g(\mathbf{x})\,d\mathbf{x}\approx\frac{\left| T\right|}{3}\sum_{i=1}^3g(\mathbf{x}_m^i)
\end{equation*}
for computing the integrals (\ref{data:b}) and (\ref{initial:condition}) involving the functions $f$ and $p_0$, respectively. Here, $\mathbf{x}_m^i$ denotes the midpoint of the $i$-th edge of $T$, for $i=1,2,3$. This formula is exact for polynomials of degree 2. Finally, the line integrals (\ref{data:a}) for the Dirichlet boundary condition $g_D$ are approximated by Simpson's rule.

\begin{table}[t]
	{\footnotesize
		\begin{center}
		\begin{tabular}{c|c|c|c|c}
		$\|r_hp-P_h\|_{\ell^{\infty}(\ell^2)} $	& $\ell=3$ &  $\ell=4$ & $\ell=5$  &  $\ell=6$ \\\hline
		$\tau=4.0$E-1 & 1.3948E-2 & 8.7120E-3 & 7.5778E-3 & 7.7539E-3\\%\hline
		$\tau=2.0$E-1 & 9.9587E-3 & 2.9405E-3 & 2.1206E-3 & 1.9511E-3\\%\hline
		$\tau=1.0$E-1 & 1.0282E-2 & 2.7437E-3 & 7.5382E-4 & 4.2576E-4\\%\hline
		$\tau=5.0$E-2 & 1.0264E-2 & 2.7232E-3 & 7.0609E-4 & 1.8993E-4\\%\hline
		$\tau=2.5$E-2 & 1.0267E-2 & 2.7211E-3 & 7.0099E-4 & 1.7894E-4\\
		\end{tabular}
		\caption{Pressure errors in $\ell^{\infty}(\ell^2)$-norm obtained for a nested collection of spatial meshes $\mathcal{T}_h^{\ell}$, with $\ell=3,4,5,6$, and various decreasing time steps $\tau$.}\label{table:ex1:pressure1}\end{center}
	}
\end{table}

\begin{table}[t]
	{\footnotesize
		\begin{center}
		\begin{tabular}{c|c|c|c|c}
		$\|r_hp-P_h\|_{\ell^{\infty}(\ell^{\infty})}$& $\ell=3$ &  $\ell=4$ & $\ell=5$  &  $\ell=6$ \\\hline
		$\tau=4.0$E-1 & 3.2717E-2 & 2.0340E-2 & 1.8256E-3 & 1.8446E-2\\%\hline
		$\tau=2.0$E-1 & 3.0892E-2 & 9.3021E-3 & 5.0352E-3 & 4.6947E-3\\%\hline
		$\tau=1.0$E-1 & 3.2036E-2 & 9.4528E-3 & 2.6495E-3 & 1.0058E-3\\%\hline
		$\tau=5.0$E-2 & 3.2013E-2 & 9.4278E-3 & 2.6247E-3 & 7.1625E-4\\%\hline
		$\tau=2.5$E-2 & 3.2083E-2 & 9.4370E-3 & 2.6212E-3 & 7.1009E-4\\
		\end{tabular}
		\caption{Pressure errors in $\ell^{\infty}(\ell^{\infty})$-norm obtained for a nested collection of spatial meshes $\mathcal{T}_h^{\ell}$, with $\ell=3,4,5,6$, and various decreasing time steps $\tau$.}\label{table:ex1:pressure2}\end{center}
	}
\end{table}

Tables \ref{table:ex1:pressure1}-\ref{table:ex1:flux2} show the pressure and velocity errors for a nested collection of spatial meshes $\mathcal{T}_h^{\ell}$, with $\ell=3,4,5,6$, and several decreasing time steps $\tau$. In particular, Tables \ref{table:ex1:pressure1} and \ref{table:ex1:pressure2} display the pressure errors computed with the norms $\|\cdot\|_{\ell^{\infty}(\ell^2)}$ and $\|\cdot\|_{\ell^{\infty}(\ell^{\infty})}$, respectively, as defined above. In both cases, we can observe an unconditionally stable behaviour of the algorithm, which converges irrespective of the size of the parameters $h$ and $\tau$ under consideration. Furthermore, as revealed by the last row in both tables, the ratios of subsequent errors imply second-order convergence in space. Accordingly, the first errors in the last column show second-order convergence in time. Note that the last two positions in this column are meaningless, since they virtually represent the error due to the spatial discretization. On the other hand, Table \ref{table:ex1:flux1} displays the computed errors for the post-processed fluxes according to formula (\ref{norm:postproc:flux}). In this case, the method shows an unconditionally convergent behaviour, with second-order convergence in time and close to second-order convergence in space. Finally, in Table \ref{table:ex1:flux2}, we observe the errors for the normal fluxes given by (\ref{norm:normal:flux}). Once again, the method is unconditionally convergent; specifically, the order of convergence approaches an asymptotic value of 2 in time and approximately $1.6$ in space, in accordance with the numerical results reported in \cite{dup:kee:98}.

\begin{table}[t]
	{\footnotesize
		\begin{center}
		\begin{tabular}{c|c|c|c|c}
		$\|\mathbf{u}-\mathbf{R}(\mathbf{u}_h)\|_{\ell^{\infty}(L^2)}$	& $\ell=3$ &  $\ell=4$ & $\ell=5$  &  $\ell=6$ \\\hline
		$\tau=4.0$E-1 & 5.6553E-1 & 2.4799E-1 & 1.6123E-1 & 1.5100E-1\\%\hline
		$\tau=2.0$E-1 & 4.7129E-1 & 1.4528E-1 & 5.5718E-2 & 4.0916E-2\\%\hline
		$\tau=1.0$E-1 & 4.9077E-1 & 1.4819E-1 & 4.4715E-2 & 1.4333E-2\\%\hline
		$\tau=5.0$E-2 & 4.8997E-1 & 1.4771E-1 & 4.4354E-2 & 1.3802E-2\\%\hline
		$\tau=2.5$E-2 & 4.8991E-1 & 1.4766E-1 & 4.4307E-2 & 1.3763E-2\\
		\end{tabular}
		\caption{Errors for the post-processed fluxes for a nested collection of spatial meshes $\mathcal{T}_h^{\ell}$, with $\ell=3,4,5,6$, and various decreasing time steps $\tau$.}\label{table:ex1:flux1}\end{center}
	}
\end{table}

\begin{table}[t]
	{\footnotesize
		\begin{center}
		\begin{tabular}{c|c|c|c|c}
		$\|\left(\mathbf{u}-\mathbf{u}_h\right)\cdot\mathbf{n}\|_{\ell^{\infty}(\ell^2)}$	& $\ell=3$ &  $\ell=4$ & $\ell=5$  &  $\ell=6$ \\\hline
		$\tau=4.0$E-1 & 1.0447E-1 & 5.8187E-2 & 5.2154E-2 & 5.2289E-2\\%\hline
		$\tau=2.0$E-1 & 8.5800E-2 & 2.9946E-2 & 1.4575E-2 & 1.2743E-2\\%\hline
		$\tau=1.0$E-1 & 8.9694E-2 & 3.0881E-2 & 1.0513E-2 & 3.6249E-3\\%\hline
		$\tau=5.0$E-2 & 8.9650E-2 & 3.0849E-2 & 1.0482E-2 & 3.5917E-3\\%\hline
		$\tau=2.5$E-2 & 8.9639E-2 & 3.0841E-2 & 1.0476E-2 & 3.5872E-3\\
		\end{tabular}
		\caption{Errors for the normal fluxes obtained for a nested collection of spatial meshes $\mathcal{T}_h^{\ell}$, with $\ell=3,4,5,6$, and various decreasing time steps $\tau$.}\label{table:ex1:flux2}\end{center}
	}
\end{table}

\subsubsection{The case of discontinuous coefficients}

\begin{figure}[t]
	\begin{center}%\hspace*{-1.8cm}
		\includegraphics[scale=0.18]{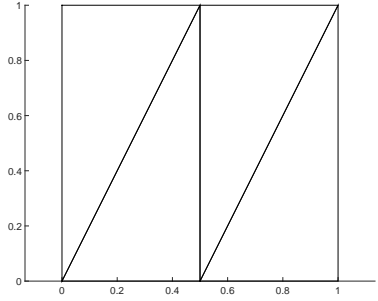}\end{center}
	\vspace*{-0.5cm}
	\caption{Coarse grid $\mathcal{T}_H$ for the Mackinnon and Carey test example with discontinuous permeability tensor.}\label{figure:mesh:mackinnon:carey}
\end{figure}

As a second example, we consider a test problem inspired by a model by Mackinnon and Carey (cf. \cite{mac:car:88}) that involves a discontinuous permeability tensor. In particular, let (\ref{ibvp}) be posed on the unit square, with $t_f=3$ and $\Gamma_D=\partial \Omega$. In this case, $K(x,y)=k(x,y)\,I$, where $I$ is the $2\times2$ identity matrix and $k(x,y)=k_1=1$, for $x\leq\frac{1}{2}$, and $k(x,y)=k_2=2$, for $x>\frac{1}{2}$. The data functions $f$, $p_0$ and $g_D$ are defined in such a way that the piecewise quadratic function
\begin{equation}\label{exact:mc}
p(x,y,t)=
\begin{cases}
\,t^2\left(a_1\dfrac{x^2}{2}+b_1x\right),&\hbox{if}\,\,0\leq x\leq\dfrac{1}{2},\\[2ex]
\,t^2\left(a_2\dfrac{x^2}{2}+b_2x+c_2\right),&\hbox{if}\,\,\dfrac{1}{2}< x\leq 1,
\end{cases}
\end{equation}
is the exact solution of the problem, where
$$
a_1=-\frac{1}{k_1}, \  a_2=-\frac{1}{k_2},\  b_1=-\frac{3a_2+a_1}{4}\frac{k_2}{k_1+k_2},\ b_2=\frac{k_1}{k_2}\,b_1,\  c_2=-b_2-\frac{a_2}{2}.
$$
The domain $\Omega$ is decomposed into $s=4$ triangular subdomains which are aligned with the permeability discontinuity, as shown in Figure \ref{figure:mesh:mackinnon:carey}. The integrals arising in the computations are approximated by the same quadrature rules as in the previous example.

Table \ref{table:mackinnon:carey} shows the pressure and velocity errors obtained for a nested collection of spatial meshes $\mathcal{T}_h^{\ell}$, with $\ell=1,2,3,4,5$, and a time step $\tau=1$E-01, considering the norms defined above. The use of the Crank--Nicolson time integrator for an exact solution of the form (\ref{exact:mc}) permits us to obtain negligible errors in time and to study the spatial convergence of the method. The last row in the table contains the average orders of convergence. As expected, the method is second-order convergent for pressures and close to second-order convergent for post-processed and normal fluxes.

\begin{table}[t]
	{\footnotesize
\begin{tabular}{c|c|c|c|c}
	&$\|r_hp-P_h\|_{\ell^{\infty}(\ell^2)} $& $\|r_hp-P_h\|_{\ell^{\infty}(\ell^{\infty})}$&$\|\mathbf{u}-\mathbf{R}(\mathbf{u}_h)\|_{\ell^{\infty}(L^2)}$ & $\|\left(\mathbf{u}-\mathbf{u}_h\right)\cdot\mathbf{n}\|_{\ell^{\infty}(\ell^2)}$\\\hline
	$\ell=1$ & 2.4218E-2 & 7.2240E-2 & 3.7393E-1 & 2.7329E-1\\%\hline
	$\ell=2$ & 5.9113E-3 & 1.8738E-2 & 1.0901E-1 & 8.0644E-2\\%\hline
	$\ell=3$ & 1.4637E-3 & 4.6708E-3 & 3.1850E-2 & 2.3434E-2\\%\hline
	$\ell=4$ & 3.6745E-4 & 1.1656E-3 & 9.5547E-3 & 6.7958E-3\\%\hline
	$\ell=5$ & 9.2331E-5 & 2.9124E-4 & 2.9657E-3 & 2.0113E-3\\\hline
	order & 2.009 & 1.989 & 1.745 & 1.772 \\
\end{tabular}
\caption{Pressure and velocity errors for the Mackinnon and Carey test example. A nested collection of spatial meshes $\mathcal{T}_h^{\ell}$, with $\ell=1,2,3,4,5$, and a time step $\tau=1$E-01 are considered. The last row in the table contains the average orders of convergence.}\label{table:mackinnon:carey}
}
\end{table}

\subsection{Numerical examples of time-dependent Darcy flow}

\subsubsection{A flow domain with low-permeability regions}

In this example, we consider (\ref{ibvp}) posed on the unit square, with $t_f=5$, $f(x,y,t)=0$ and $p_0(x,y)=1-x$. Furthermore, $\Gamma_D=\{0,1\}\times(0,1)$ and $\Gamma_N=(0,1)\times\{0,1\}$. The pressure is specified to be equal to 1 on the left boundary and equal to 0 on the right boundary. A zero-flux condition is imposed on $\Gamma_N$. The flow domain contains two low-permeability regions, namely, $R_1=(0.2,0.3)\times(0,0.8)$ and $R_2=(0.6,0.7)\times(0.3,1)$. In particular, the permeability tensor is $K=k(x,y)\,I$, where $I$ is the $2\times2$ identity matrix and
$$k(x,y)=
\begin{cases}
\,1,& \hbox{if}\,\,(x,y)\in \Omega\setminus(R_1\cup R_2),\\
\,10^{-6},& \hbox{if}\,\,(x,y)\in R_1\cup R_2.
\end{cases}
$$
This numerical example permits us to test the behaviour of the algorithm on problems involving abrupt variations in the permeability. A stationary version of this problem was considered in \cite{gan:lis:sha:zen:02} in the framework of heat flow applications.

\begin{figure}[t]
	\begin{center}\hspace*{-1.9cm}
		\begin{minipage}[t]{0.48\textwidth}
			\begin{center}\includegraphics[scale=0.18]{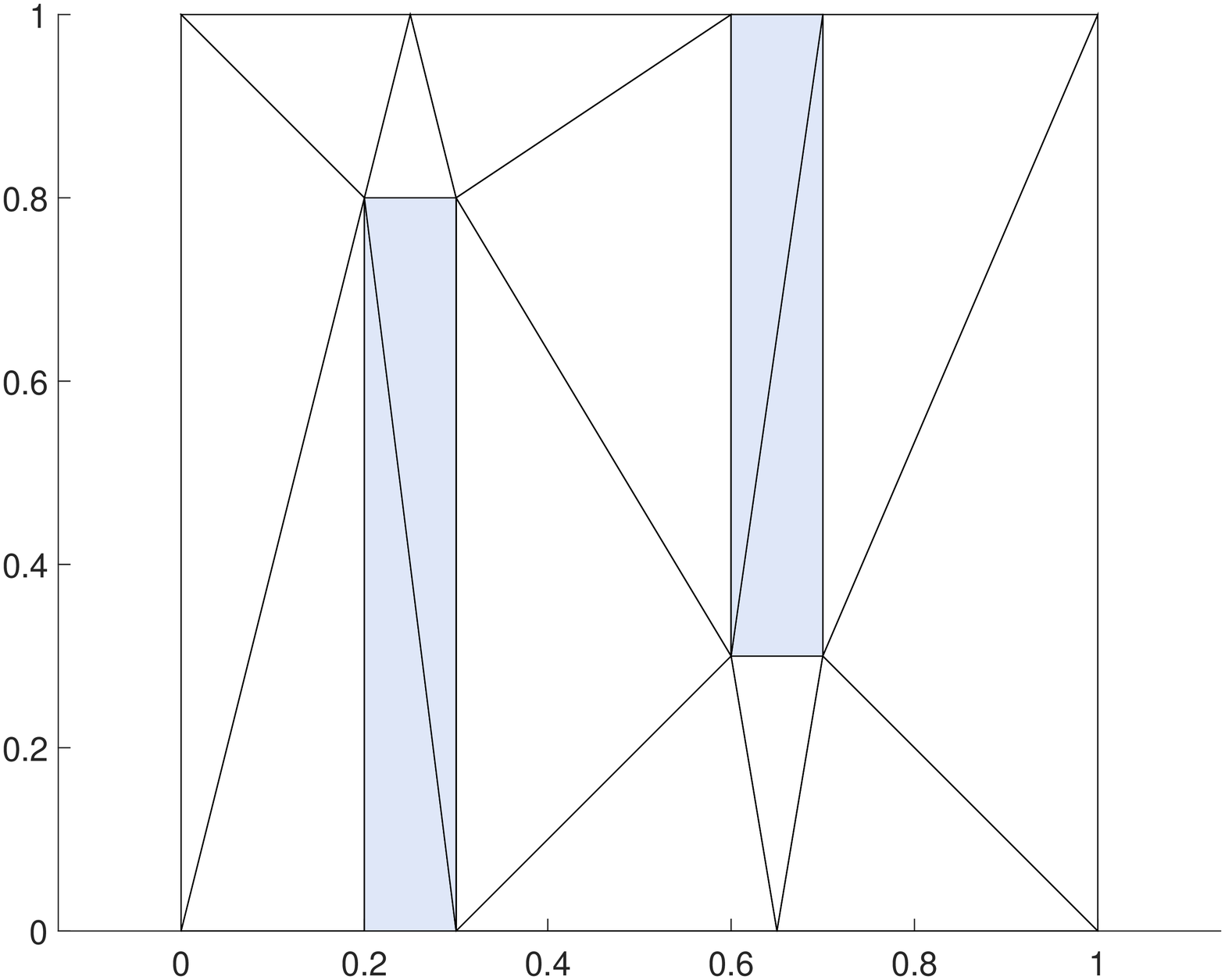}\end{center}
		\end{minipage}
		\hspace*{-0.15cm}
		\begin{minipage}[t]{0.48\textwidth}\begin{center}\includegraphics[scale=0.18]{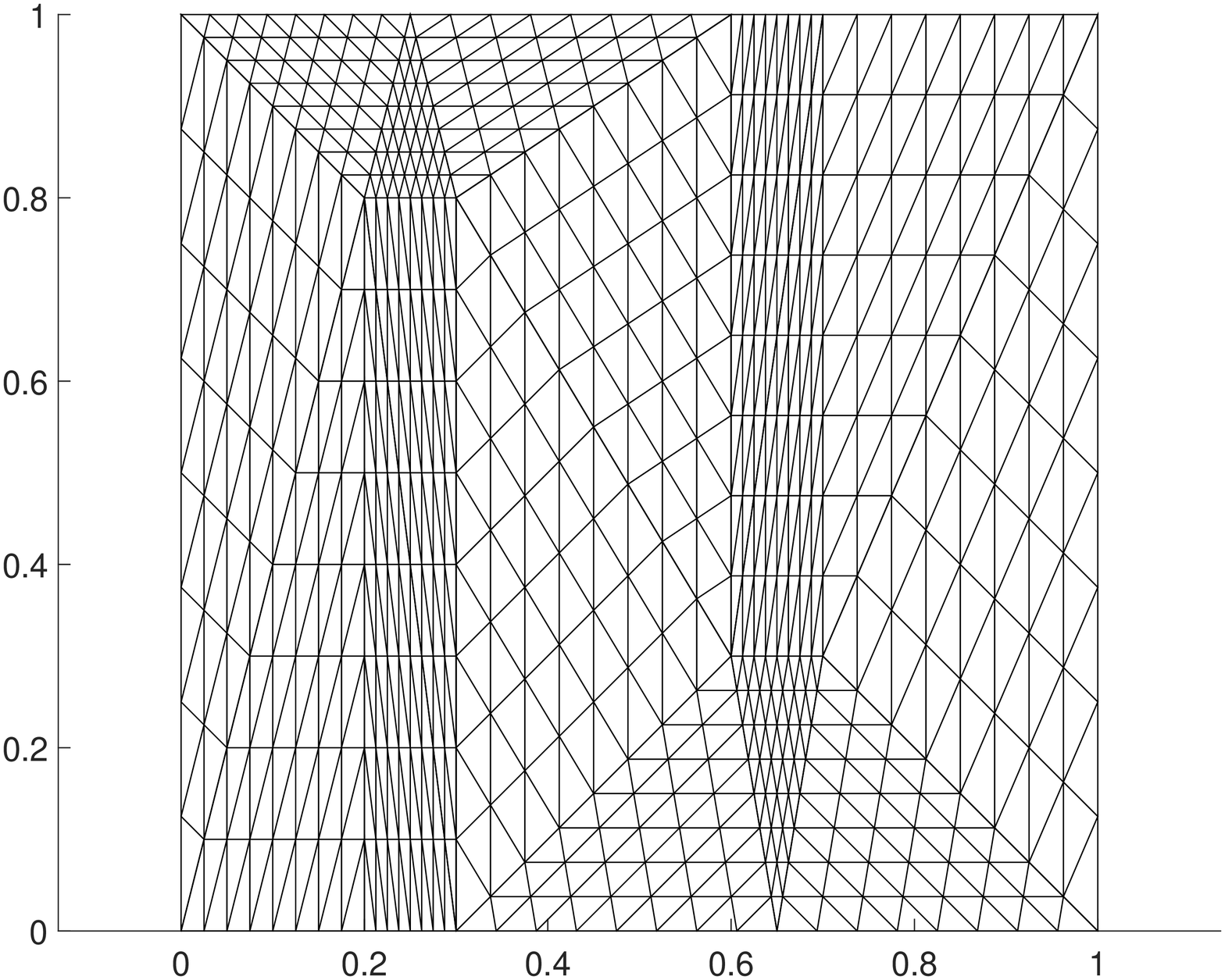}\end{center}
		\end{minipage}
	\end{center}
	\vspace*{-0.5cm}
	\caption{Coarse grid $\mathcal{T}_H$ (left) and fine grid $\mathcal{T}_h\equiv\mathcal{T}_h^3$ (right) for the unit square with two low-permeability regions.}\label{figure:meshes:two:fingers}
\end{figure}

Figure~\ref{figure:meshes:two:fingers} (left) shows the coarse grid $\mathcal{T}_H$ used in this example. It consists of $s=16$ triangular subdomains and it is adapted to the geometry of the low-permeability regions, depicted in the plot by shaded areas. In turn, Figure \ref{figure:meshes:two:fingers} (right) shows the fine grid $\mathcal{T}_h\equiv\mathcal{T}_h^3$ obtained after three regular refinement processes.

\begin{figure}[t]
	\begin{center}\hspace*{-1.9cm}
		\begin{minipage}[t]{0.48\textwidth}
            \begin{center}\includegraphics[scale=0.18]{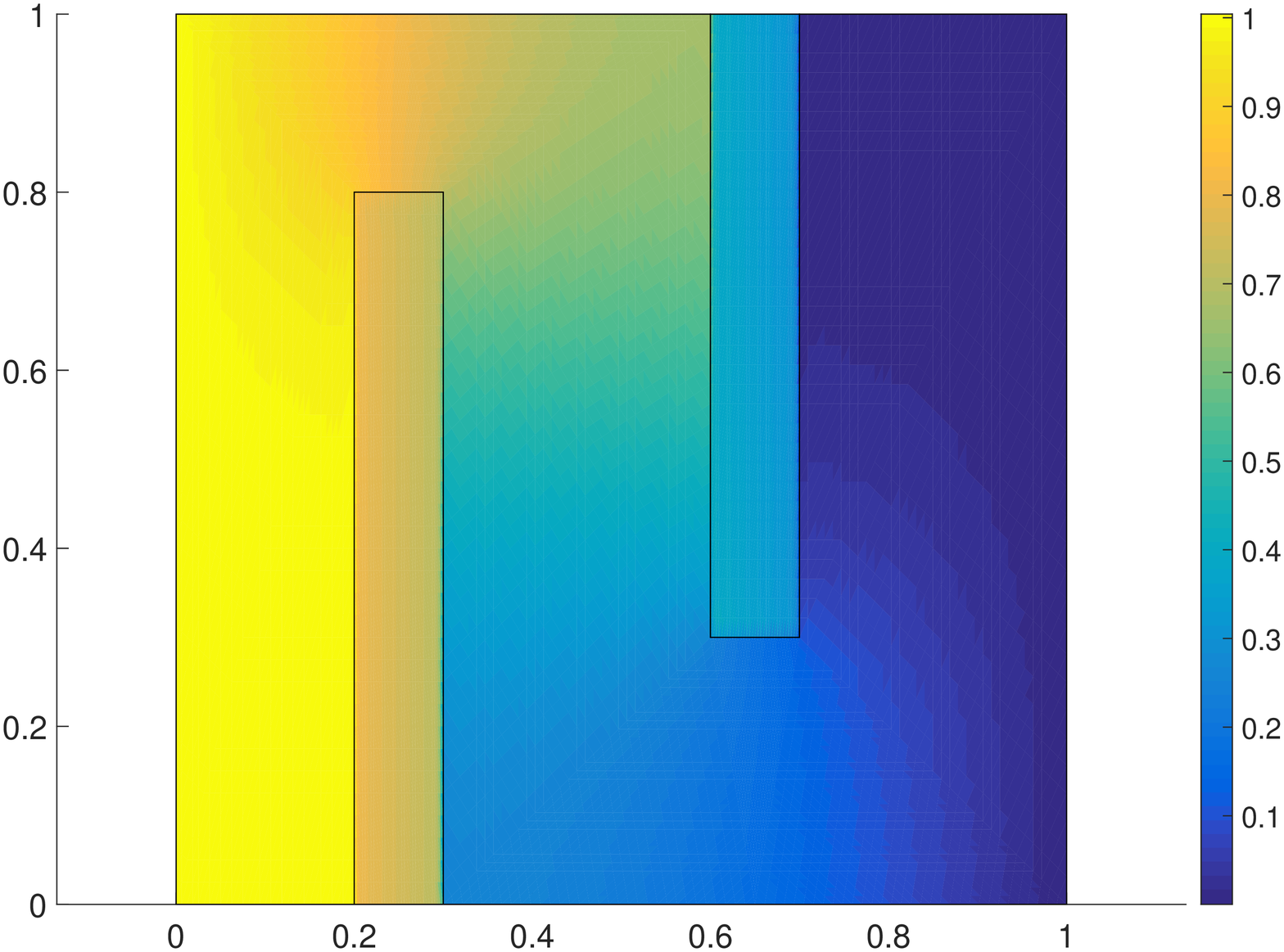}\end{center}
		\end{minipage}
		\hspace*{-0.15cm}
		\begin{minipage}[t]{0.48\textwidth}\begin{center}\includegraphics[scale=0.18]{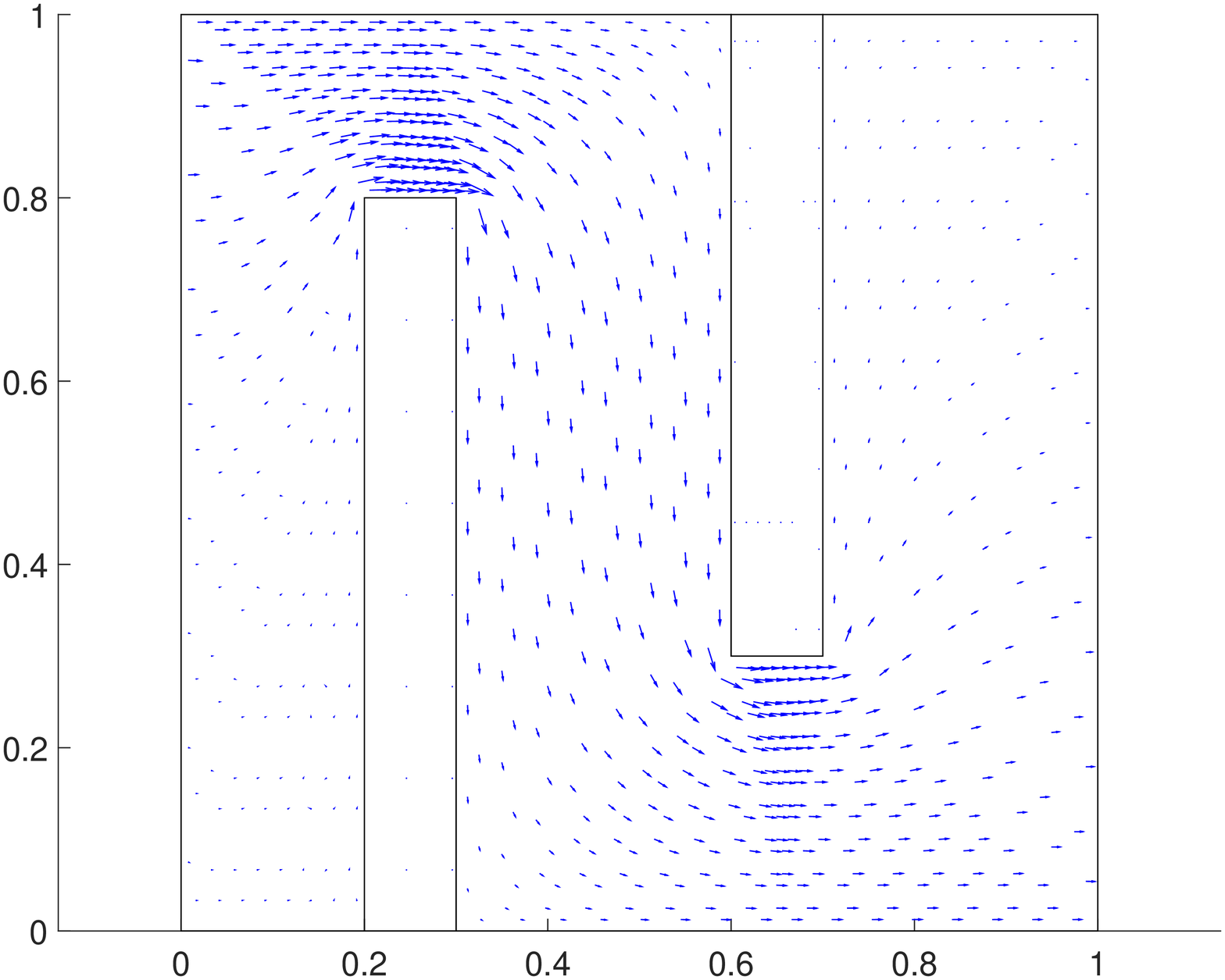}\end{center}
		\end{minipage}
	\end{center}
	\vspace*{-0.5cm}
	\caption{Pressure distribution (left) and velocity field (right) for the Example 5.2.1.}\label{figure:pressure:flux:two:fingers}
\end{figure}

In Figure~\ref{figure:pressure:flux:two:fingers} (left), we display the pressure distribution obtained  once the stationary state is reached. In this case, we consider the triangular mesh $\mathcal{T}_h^5$, obtained after five regular refinement processes of the coarse grid $\mathcal{T}_H$, and a time step $\tau=5$E-02. This figure illustrates the effect of the low-permeability regions on the pressure distribution. On the other hand, Figure~\ref{figure:pressure:flux:two:fingers} (right) shows the velocity field obtained at the stationary state. As usual, the length of the arrows is proportional to the module of the vectors. In this case, we consider the triangular mesh $\mathcal{T}_h^3$ shown in Figure \ref{figure:meshes:two:fingers} (right) and a time step $\tau=5$E-02. As expected from the physical configuration, no flow enters the low-permeability regions.

\subsubsection{A flow domain with holes}

\begin{figure}[t]
	\begin{center}\hspace*{-1.9cm}
		\begin{minipage}[t]{0.48\textwidth}
			\begin{center}\includegraphics[scale=0.18]{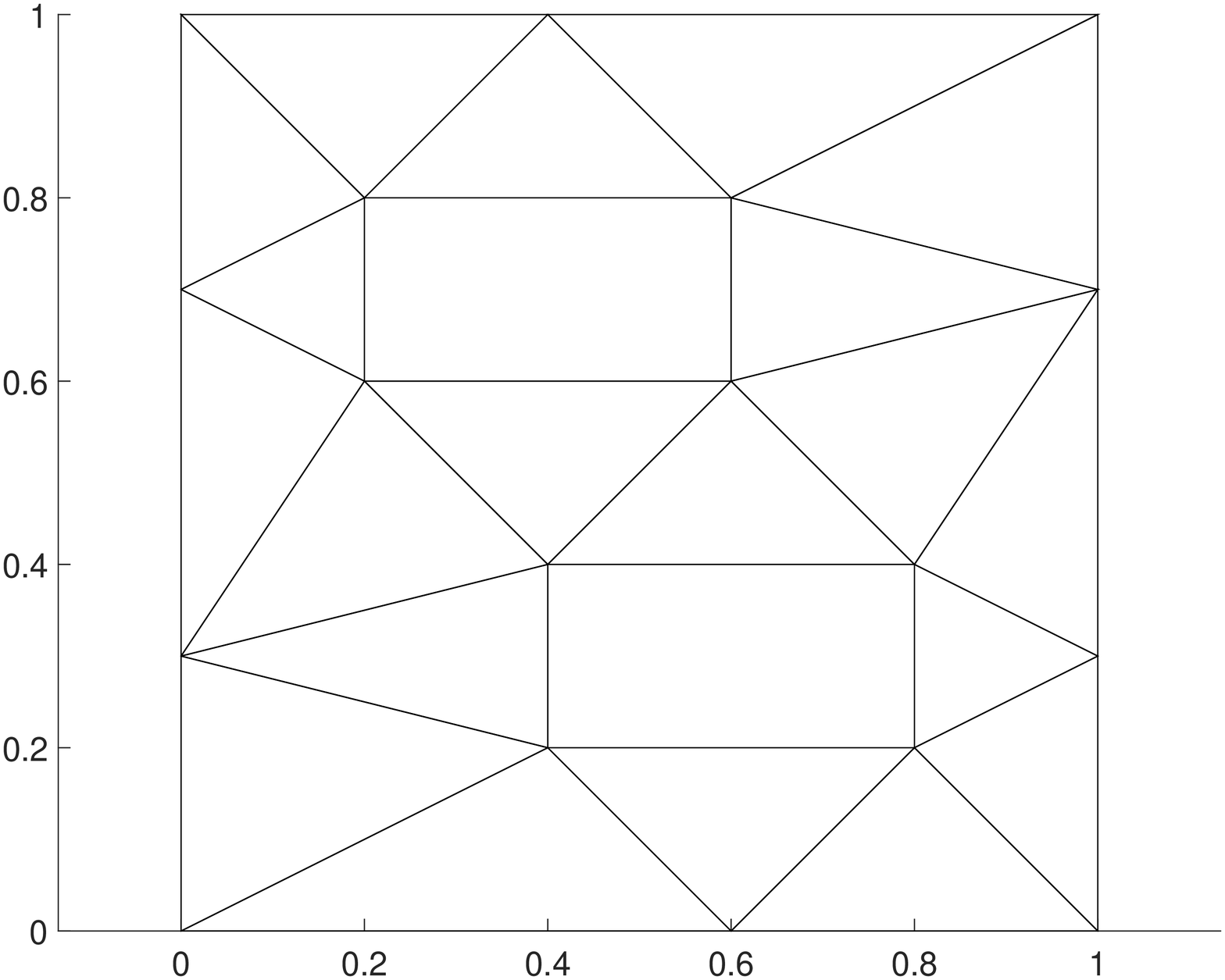}\end{center}
		\end{minipage}
		\hspace*{-0.15cm}
		\begin{minipage}[t]{0.48\textwidth}\begin{center}\includegraphics[scale=0.18]{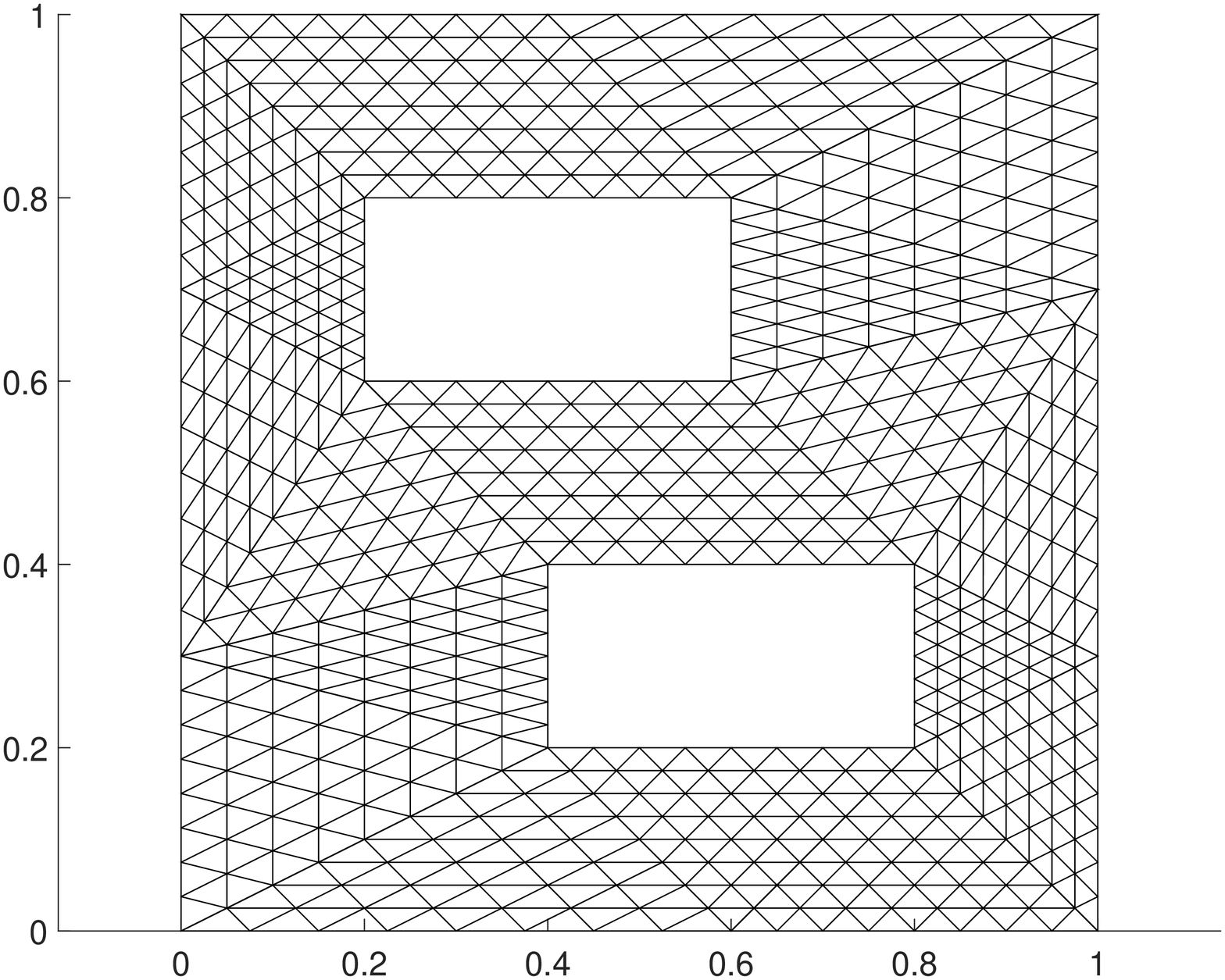}\end{center}
		\end{minipage}
	\end{center}
	\vspace*{-0.5cm}
	\caption{Coarse grid $\mathcal{T}_H$ (left) and fine grid $\mathcal{T}_h\equiv\mathcal{T}_h^3$ (right) for the non-simply connected domain $\Omega$.}\label{figure:meshes:two:holes}
\end{figure}

\begin{figure}[t]
	\begin{center}\hspace*{-1.9cm}
		\begin{minipage}[t]{0.48\textwidth}
			\begin{center}\includegraphics[scale=0.18]{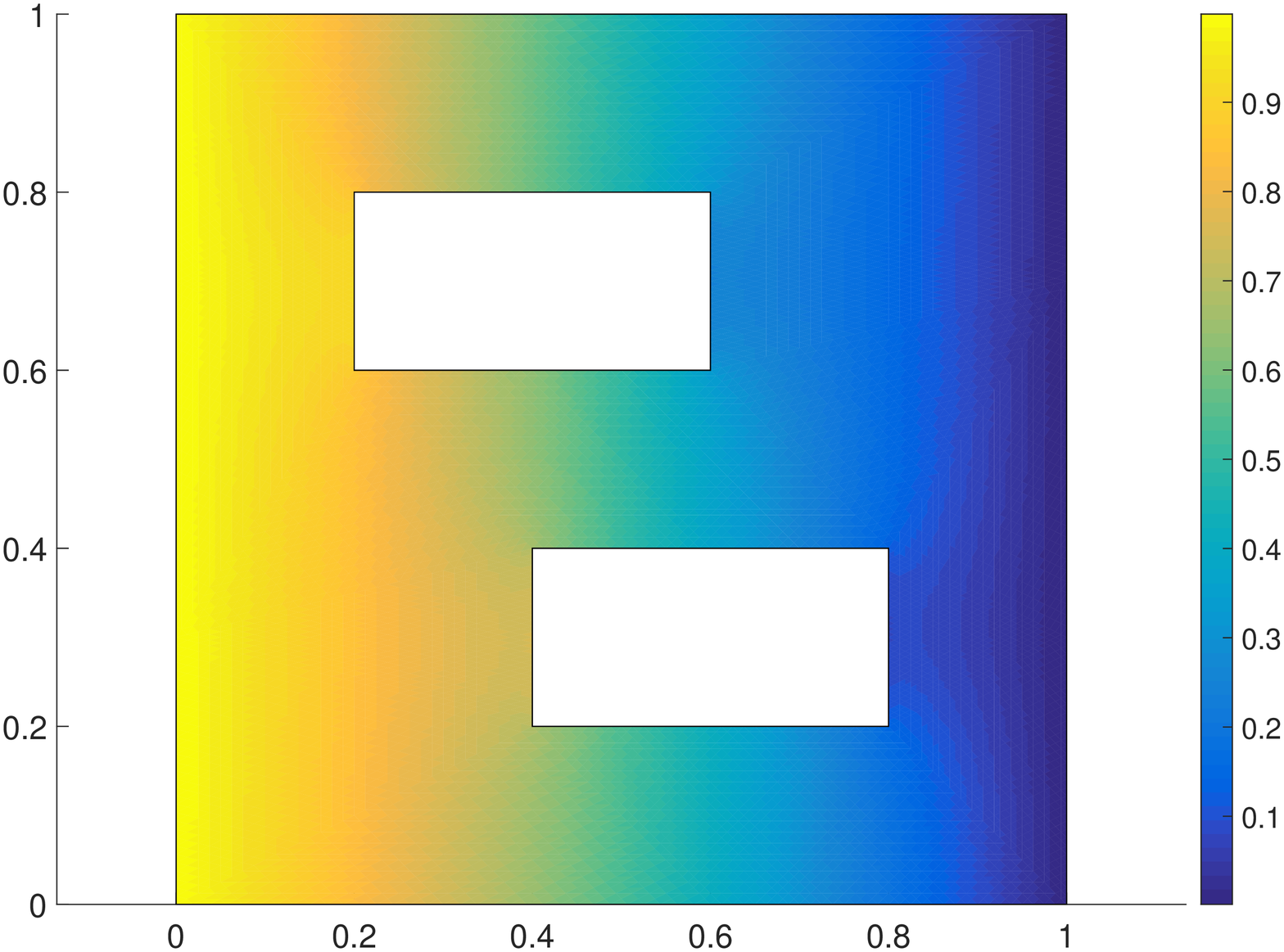}\end{center}
		\end{minipage}
		\hspace*{-0.15cm}
		\begin{minipage}[t]{0.48\textwidth}\begin{center}\includegraphics[scale=0.18]{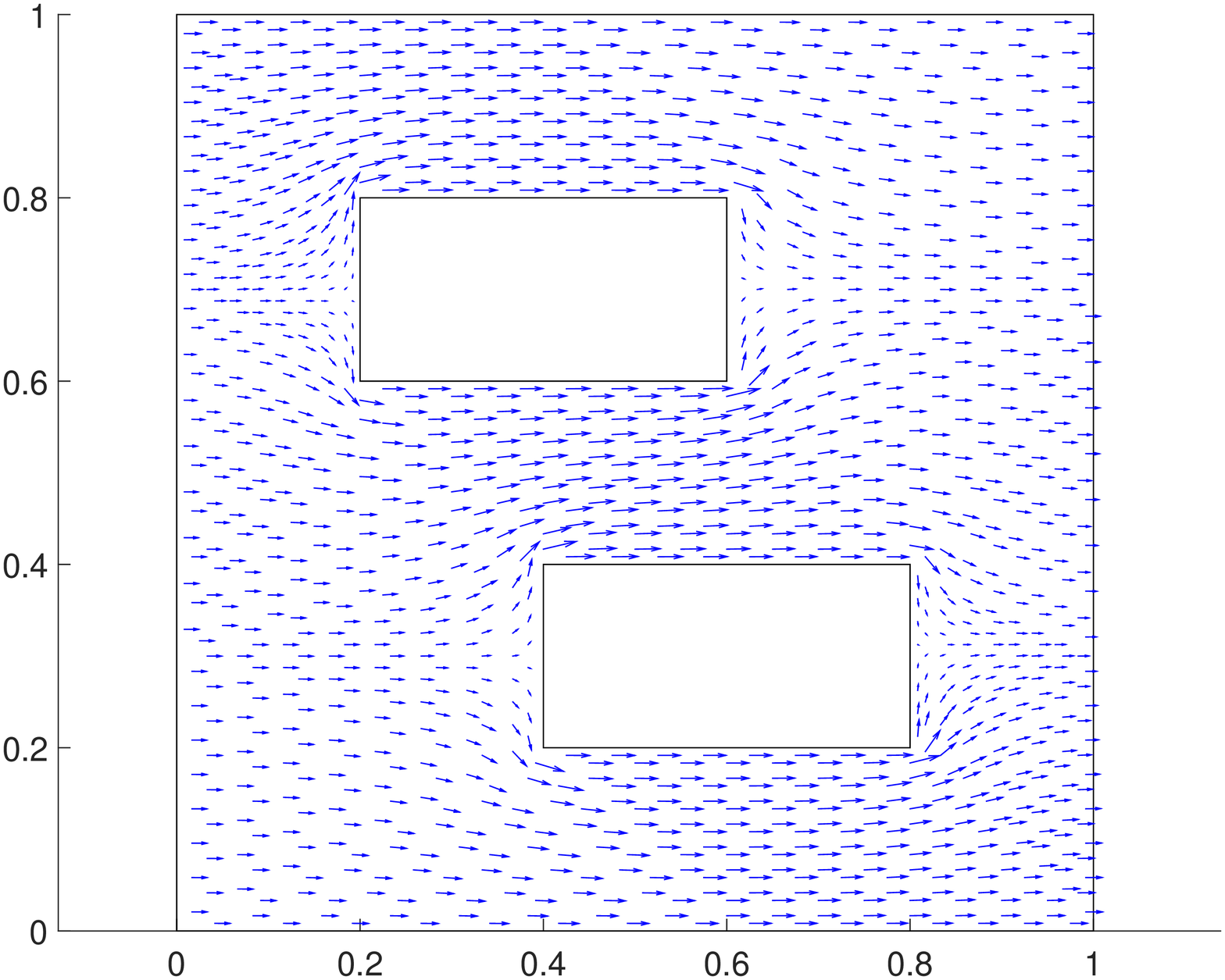}\end{center}
		\end{minipage}
	\end{center}
	\vspace*{-0.5cm}
	\caption{Pressure distribution (left) and velocity field (right) for the Example 5.2.2.}\label{figure:pressure:flux:two:holes}
\end{figure}

Finally, we consider (\ref{ibvp}) posed on $\Omega\times(0,t_f]$, where $\Omega$ is a non-simply connected domain with two rectangular holes, and $t_f=5$. Specifically, we define $\Omega=(0,1)^2\setminus(R_1\cup R_2)$, where $R_1=[0.2,0.6]\times[0.6,0.8]$ and $R_2=[0.4,0.8]\times[0.2,0.4]$. The source/sink term is $f(x,y,t)=0$, and the initial condition is $p_0(x,y)=1-x$. The boundary of $\Omega$ is decomposed into $\Gamma_D=\{0,1\}\times(0,1)$ and $\Gamma_N=(0,1)\times\{0,1\}\cup\partial R_1\cup\partial R_2$. In this context, the pressure is set to be equal to 1 on the boundary $\{0\}\times(0,1)$ and equal to 0 on the boundary $\{1\}\times(0,1)$. In turn, zero flux is imposed on $\Gamma_N$. In this example, we consider a uniform and isotropic permeability tensor $K=I$, where $I$ is the $2\times2$ identity matrix. This experiment permits us to test the behaviour of the method on problems posed on non-simply connected domains. Stationary versions of this problem considering various configurations of holes in a square domain can be found in \cite{gan:lis:sha:zen:02}.

Figure \ref{figure:meshes:two:holes} (left) shows the coarse grid $\mathcal{T}_H$ used in this example. It consists of $s=20$ triangular subdomains, which are aligned with the two rectangular holes. Figure \ref{figure:meshes:two:holes} (right) shows the fine grid $\mathcal{T}_h\equiv\mathcal{T}_h^3$ obtained after three regular refinement processes.

The computed pressures and fluxes are displayed on Figure \ref{figure:pressure:flux:two:holes}. The left-hand plot shows the pressure distribution at the stationary state, when considering the triangular mesh $\mathcal{T}_h^5$ and a time step $\tau=5$E-02. The effect of both holes on the pressure distribution is quite significant. In turn, the right-hand plot shows the velocity field obtained at the stationary state, when considering the triangular mesh $\mathcal{T}_h^3$ depicted in Figure~\ref{figure:meshes:two:holes} (right) and $\tau=5$E-02. In this case, the fluid flows around the internal boundaries of the domain.

\section{Concluding remarks}\label{sec:conclusions}

A novel algorithm for solving time-dependent Darcy problems with possibly discontinuous tensor coefficients has been designed and tested on different geometries. The method is conceived to fit in the framework of hierarchical grids, taking advantage of their flexibility to handle complex problems while preserving the efficiency of stencil-based operations.

A relevant feature of the numerical procedure is that, provided that the material properties are constant within each subdomain, the corresponding stencil remains unchanged regardless of the level of refinement. In addition, the fully discrete scheme yields a local problem on each subdomain whose solution does not depend on the neighbouring subdomains, thus allowing for parallelization. Since the subdomain problems are all the same size, a perfectly balanced workload can be assigned among parallel processors in the solution process.

The structure of the resulting linear systems can lead to the development of very efficient linear solvers based on a judicious combination of preconditioning techniques and multigrid methods. Furthermore, the ease of implementation and robust numerical behaviour of the algorithm can make it suitable for more general flow problems. Both topics are the scope of current research.

\appendix

\section{Derivation of the stencil for $MP_h$}\label{sec:appendix}

In this appendix, we describe how to deduce the coefficients of the local 10-point stencil associated to $MP_h$, when $K$ is a tensor whose components do not depend on the spatial variables. Since matrix $M$ decouples across subdomains (cf. Subsection \ref{subsec:ccfd}), it suffices to study each subdomain separately.

\begin{figure}[t]
	\begin{center}
		\hspace*{8.5ex}
		\unitlength=0.06cm
		\begin{picture}(260,110)
		%\allinethickness{0.75pt}
		%Ponemos el nombre del subdominio
		\put(172,65){{\small$\Omega_k$}}
		%Dibujamos los ejes de coordenadas
		\put(22,0){\vector(1,0){173}}
		\put(32,-5){\vector(-1,2){56}}
		%Dibujamos el mallado
		%L\~{A}要eas horizontales
		
		\put(50,20){\color{brown}\line(1,0){120}}
		\put(70,40){\line(1,0){30}}\put(130,40){\line(1,0){30}}
		\put(90,60){\line(1,0){60}}
		\put(110,80){\line(1,0){30}}
		\put(100,40){\color{brown}\line(1,0){30}}
		
		%L\~{A}要eas con pendiente 1
		\put(50,20){\color{red}\line(1,1){80}}
		\put(80,20){\line(1,1){20}}\put(120,60){\line(1,1){20}}
		\put(110,20){\line(1,1){40}}
		\put(140,20){\line(1,1){20}}
		\put(100,40){\color{red}\line(1,1){20}}

		%L\~{A}要eas con pendiente negativa
		\put(80,20){\line(-1,2){10}}
		\put(110,20){\line(-1,2){10}}
		\put(140,20){\line(-1,2){30}}
		\put(170,20){\color{blue}\line(-1,2){40}}
		\put(100,40){\color{blue}\line(-1,2){10}}
		
		%Ponemos nombres a los ejes
		\put(194,3){{\footnotesize$\mathbf{e}_1'$}}
		\put(-20,107){{\footnotesize$\mathbf{e}_2'$}}
		
		%Dibujamos los ticks en los ejes
%				\Dotline(60,0)(50,20){1.5}
%				\Dotline(90,0)(80,20){1.5}
%				\Dotline(120,0)(110,20){1.5}
%				\Dotline(150,0)(140,20){1.5}
%				\Dotline(180,0)(170,20){1.5}
%				%
%				\Dotline(20,20)(50,20){1.5}
%				\Dotline(10,40)(70,40){1.5}
%				\Dotline(0,60)(90,60){1.5}
%				\Dotline(-10,80)(110,80){1.5}
%				\Dotline(-20,100)(130,100){1.5}
		
		%Ponemos nombre a los ticks de los ejes
		\put(61,-6){{\footnotesize$1$}}
		\put(91,-6){{\footnotesize$2$}}
		\put(121,-6){{\footnotesize$3$}}
		\put(151,-6){{\footnotesize$4$}}
		\put(181,-6){{\footnotesize$5$}}
		\put(13,18){{\footnotesize$1$}}
		\put(3,38){{\footnotesize$2$}}
		\put(-6.5,58){{\footnotesize$3$}}
		\put(-16.5,78){{\footnotesize$4$}}
		\put(-26,98){{\footnotesize$5$}}
		
		%Ponemos nombres a los nodos de la malla
		\put(38,21){{\scriptsize$(1,1)$}}
		\put(76,14){{\scriptsize$(2,1)$}}
		\put(106,14){{\scriptsize$(3,1)$}}
		\put(136,14){{\scriptsize$(4,1)$}}
		\put(171,21){{\scriptsize$(5,1)$}}
		\put(58,41){{\scriptsize$(2,2)$}}
		\put(85.5,42){{\scriptsize$(3,2)$}}
		\put(98.4,38){$\bullet$}
		\put(136,42){{\scriptsize$(4,2)$}}
		\put(78,61){{\scriptsize$(3,3)$}}
		\put(105.5,62){{\scriptsize$(4,3)$}}
		\put(98,81){{\scriptsize$(4,4)$}}
		\put(161,41){{\scriptsize$(5,2)$}}
		\put(151,61){{\scriptsize$(5,3)$}}
		\put(141,81){{\scriptsize$(5,4)$}}
		\put(125,103){{\scriptsize$(5,5)$}}
		
		%		%Ponemos nombres a los lados de la malla
		\put(111,47){{\color{red}\footnotesize$e_{3,2}^2$}}
		\put(111,33.5){{\color{brown}\footnotesize$e_{3,2}^3$}}
		\put(95.5,51){{\color{blue}\footnotesize$e_{3,2}^1$}}
				
		\end{picture}
	\end{center}\vspace*{-0.3cm}
	\caption{A mesh $\mathcal{T}_{h,k}^2$ covering the subdomain $\Omega_k$ is represented in an oblique coordinate system $\{\mathbf{e}_1',\mathbf{e}_2'\}$. The local numbering of grid points is included. In addition, the three edges $e_{3,2}^1$, $e_{3,2}^2$ and $e_{3,2}^3$ associated to the grid point $(3,2)$ are depicted in blue, red and brown, respectively.}
	\label{fig:new:ref:sys}
\end{figure}
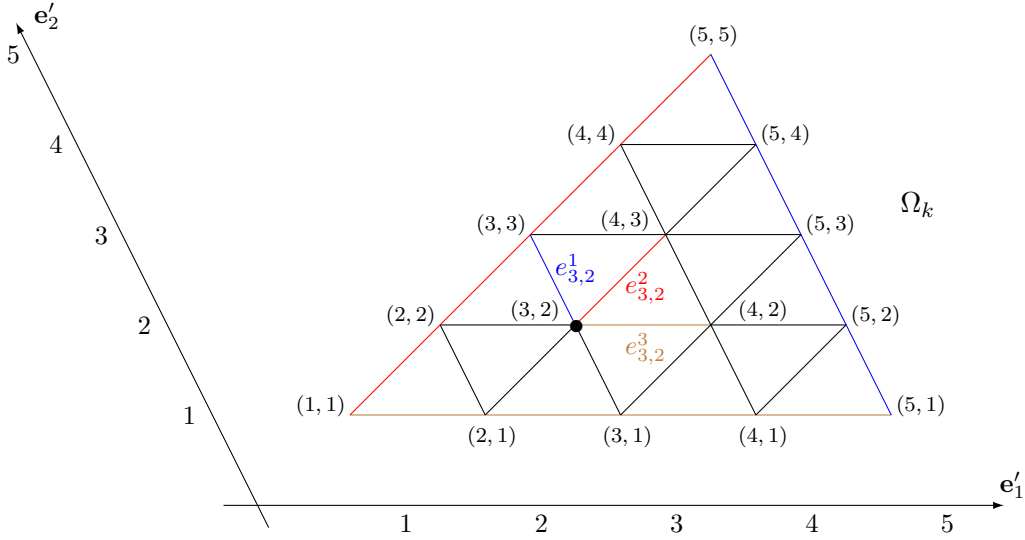

\begin{figure}[t]
	\begin{center}
		%		\hspace*{8.5ex}
		\unitlength=0.08cm
		\begin{picture}(160,33)
		%\allinethickness{0.75pt}
		
		%% Figura de la izquierda
		
		%Pintamos el fondo del triangulo up
				{\color{brown!70}\polygon*(26.75,5)(56.75,5)(46.75,25)}
		%Pintamos el fondo del triangulo down
				{\color{blue!15}\polygon*(25,5)(45,25)(15,25)}
		%Dibujamos los dos tri\'{a}ngulos
		%L\~{A}要eas horizontales
		\put(25,5){\line(1,0){30}}
		\put(15,25){\line(1,0){30}}
		%L\~{A}要eas con pendiente 1
		\put(25,5){\line(1,1){20}}
		
		%L\~{A}要eas con pendiente negativa
		\put(25,5){\line(-1,2){10}}
		\put(55,5){\line(-1,2){10}}
		
		%Ponemos nombres a los tri\'{a}gulos
		\put(40,9){{\footnotesize$T_{i,j}^U$}}
		\put(25,18){{\footnotesize$T_{i,j}^D$}}
		
		%Ponemos nombres a los nodos de la malla
		\put(17,1){{\scriptsize$(i,j)$}}
		\put(56,1){{\scriptsize$(i+1,j)$}}
		\put(45,27){{\scriptsize$(i+1,j+1)$}}
		\put(0,26){{\scriptsize$(i,j+1)$}}
		
		%Ponemos nombres a los lados de la malla
		\put(13,15){{\footnotesize$e_{i,j}^1$}}
		\put(33,15){{\footnotesize$e_{i,j}^2$}}
		\put(39,0){{\footnotesize$e_{i,j}^3$}}
		\put(51,15){{\footnotesize$e_{i+1,j}^1$}}
		\put(27,27){{\footnotesize$e_{i,j+1}^3$}}

		%% Figura de la derecha
		
		%Pintamos el fondo del triangulo up
				{\color{brown!70}\polygon*(116.75,5)(146.75,5)(136.75,25)}
		%Pintamos el fondo del triangulo down
				{\color{blue!15}\polygon*(115,5)(135,25)(105,25)}
		%Dibujamos los dos tri\'{a}ngulos
		%L\~{A}要eas horizontales
		\put(115,5){\line(1,0){30}}
		\put(105,25){\line(1,0){30}}
		%L\~{A}要eas con pendiente 1
		\put(115,5){\line(1,1){20}}
		
		%L\~{A}要eas con pendiente negativa
		\put(115,5){\line(-1,2){10}}
		\put(145,5){\line(-1,2){10}}
		
		%Ponemos nombres a los tri\'{a}gulos
		\put(130,9){{\footnotesize$T_{i,j}^U$}}
		\put(115,18){{\footnotesize$T_{i,j}^D$}}
		
		%Numeramos nodos del tri\'{a}ngulo up en rojo
		\put(115,0){{\color{brown}\footnotesize$\mathbf{r}_1$}}
		\put(142,0){{\color{brown}\footnotesize$\mathbf{r}_2$}}
		\put(137,23){{\color{brown}\footnotesize$\mathbf{r}_3$}}
		%Numeramos nodos del tri\'{a}ngulo down en azul
		\put(131,27){{\color{blue}\footnotesize$\mathbf{r}_1$}}
		\put(105,27){{\color{blue}\footnotesize$\mathbf{r}_2$}}
		\put(109,5){{\color{blue}\footnotesize$\mathbf{r}_3$}}
		
		\end{picture}
	\end{center}\vspace*{-0.3cm}
	\caption{Representation of two neighbouring triangles $T_{i,j}^U$ and $T_{i,j}^D$. The local notation of their sides and vertices (left) and the local ordering of their vertices (right) are included.}
	\label{fig:up:down}
\end{figure}
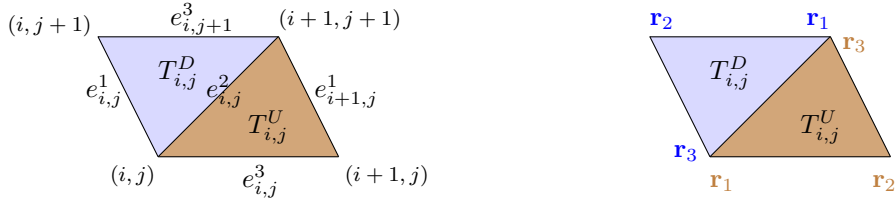

\begin{figure}[t]
	\begin{center}
		\hspace*{-3.2ex}
		\unitlength=0.06cm
		\begin{picture}(260,80)
		%\allinethickness{0.75pt}
		%Pintamos el fondo de los triangulos up
						{\color{brown!70}\polygon*(110.5,60)(140.5,60)(130.5,80)}
						{\color{brown!70}\polygon*(138,60)(168,60)(158,80)}
						{\color{brown!70}\polygon*(85.75,40)(115.75,40)(105.75,60)}
						{\color{brown!70}\polygon*(113.5,40)(143.5,40)(133.5,60)}
						{\color{brown!70}\polygon*(141.25,40)(171.25,40)(161.25,60)}
						{\color{brown!70}\polygon*(89,20)(119,20)(109,40)}
						{\color{brown!70}\polygon*(116.75,20)(146.75,20)(136.75,40)}
						
						%Pintamos el fondo de los triangulos down
						{\color{blue!15}\polygon*(94.5,60)(104.5,40)(124.5,60)}
						{\color{blue!15}\polygon*(122.25,60)(132.25,40)(152.25,60)}
						{\color{blue!15}\polygon*(100,40)(110,20)(130,40)}
		
		%Dibujamos los ejes de coordenadas
		\put(35,6){\vector(1,0){170}}
		\put(45,0){\vector(-1,2){45}}
		%Dibujamos el mallado
		%L\~{A}要eas horizontales
		\put(50,20){\line(1,0){120}}
		\put(70,40){\line(1,0){90}}
		\put(60,60){\line(1,0){120}}
		\put(80,80){\line(1,0){90}}
		%L\~{A}要eas con pendiente 1
		\put(60,60){\line(1,1){20}}
		\put(50,20){\line(1,1){60}}\put(80,20){\line(1,1){60}}
		\put(110,20){\line(1,1){60}}\put(140,20){\line(1,1){40}}
		
		%L\~{A}要eas con pendiente negativa
		\put(80,20){\line(-1,2){20}}\put(110,20){\line(-1,2){30}}
		\put(140,20){\line(-1,2){30}}\put(170,20){\line(-1,2){30}}
		\put(180,60){\line(-1,2){10}}
		
		%Ponemos nombres a los ejes
		\put(202,-2){{\footnotesize$\mathbf{e}_1'$}}
		\put(-7,85){{\footnotesize$\mathbf{e}_2'$}}
		
		%Dibujamos los ticks en los ejes
%						\Dotline(87,6)(80,20){1.5}
%						\Dotline(117,6)(110,20){1.5}
%						\Dotline(147,6)(140,20){1.5}
%						%
%						\Dotline(35,20)(50,20){1.5}
%						\Dotline(25,40)(70,40){1.5}
%						\Dotline(15,60)(60,60){1.5}
		
		%Ponemos nombre a los ticks de los ejes
		\put(83,0){{\scriptsize$i-1$}}
		\put(118,0){{\scriptsize$i$}}
		\put(143,0){{\scriptsize$i+1$}}
		\put(20,18){{\scriptsize$j-1$}}
		\put(19,38){{\scriptsize$j$}}
		\put(0,58){{\scriptsize$j+1$}}
		
		%Ponemos nombre a los tri\~{A}!`ngulos up
		\put(122,25){{\scriptsize$T_{i,j-1}^U$}}
		\put(88,24){{\scriptsize$T_{i-1,j-1}^U$}}
		\put(80,45){{\scriptsize$T_{i-1,j}^U$}}
		\put(113,45){{\scriptsize$T_{i,j}^U$}}
		\put(140,45){{\scriptsize$T_{i+1,j}^U$}}
		\put(128,64){{\scriptsize$T_{i+1,j+1}^U$}}
		\put(100,65){{\scriptsize$T_{i,j+1}^U$}}
		
		%Ponemos nombre a los tri\~{A}!`ngulos down
		\put(126.5,53){{\scriptsize$T_{i+1,j}^D$}}
		\put(99,53){{\scriptsize$T_{i,j}^D$}}
		\put(107,33){{\scriptsize$T_{i,j-1}^D$}}
		
		\end{picture}
	\end{center}\vspace*{-0.3cm}
	\caption{10-point stencil for the upward triangle $T_{i,j}^U\in\mathcal{T}_{h,k}$ represented in the oblique coordinate system $\{\mathbf{e}_1', \mathbf{e}_2'\}$.}
	\label{fig:stencil:up}
\end{figure}
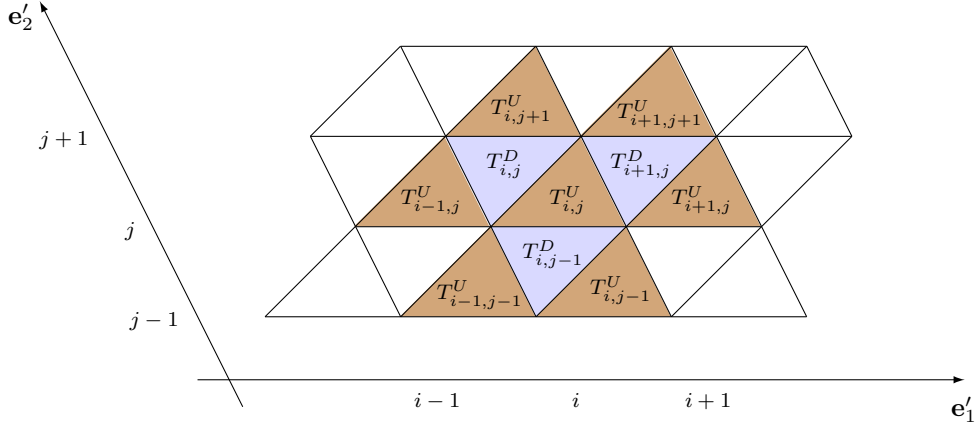

\begin{figure}[t]
	\begin{center}
		\hspace*{-3.2ex}
		\unitlength=0.06cm
		\begin{picture}(260,85)
		%\allinethickness{0.75pt}
		%Pintamos el fondo de los triangulos up
		{\color{brown!70}\polygon*(120.25,60)(150.25,60)(140.25,80)}
		{\color{brown!70}\polygon*(98,40)(128,40)(118,60)}
		{\color{brown!70}\polygon*(125.75,40)(155.75,40)(145.75,60)}
		%Pintamos el fondo de los triangulos down
		{\color{blue!15}\polygon*(103.5,80)(113.5,60)(133.5,80)}
		{\color{blue!15}\polygon*(131.25,80)(141.25,60)(161.25,80)}
		{\color{blue!15}\polygon*(79,60)(89,40)(109,60)}
		{\color{blue!15}\polygon*(106.75,60)(116.75,40)(136.75,60)}
		{\color{blue!15}\polygon*(134.5,60)(144.5,40)(164.5,60)}
		{\color{blue!15}\polygon*(82.25,40)(92.25,20)(112.25,40)}
		{\color{blue!15}\polygon*(110,40)(120,20)(140,40)}	
		%Dibujamos los ejes de coordenadas
		\put(35,6){\vector(1,0){170}}
		\put(45,0){\vector(-1,2){45}}
		%Dibujamos el mallado
		%L\~{A}要eas horizontales
		\put(60,20){\line(1,0){90}}
		\put(50,40){\line(1,0){120}}
		\put(70,60){\line(1,0){90}}
		\put(60,80){\line(1,0){120}}
		%%L\~{A}要eas con pendiente 1
		\put(50,40){\line(1,1){40}}
		\put(60,20){\line(1,1){60}}\put(90,20){\line(1,1){60}}
		\put(120,20){\line(1,1){60}}\put(150,20){\line(1,1){20}}
		
		%L\~{A}要eas con pendiente negativa
		\put(60,20){\line(-1,2){10}}\put(90,20){\line(-1,2){30}}
		\put(120,20){\line(-1,2){30}}\put(150,20){\line(-1,2){30}}
		\put(170,40){\line(-1,2){20}}
		
		%Ponemos nombres a los ejes
		\put(202,-2){{\footnotesize$\mathbf{e}_1'$}}
		\put(-7,85){{\footnotesize$\mathbf{e}_2'$}}
		
		%Dibujamos los ticks en los ejes
%				\Dotline(97,6)(90,20){1.5}
%				\Dotline(127,6)(120,20){1.5}
%				\Dotline(157,6)(150,20){1.5}
%%		%		%
%				\Dotline(35,20)(60,20){1.5}
%				\Dotline(25,40)(50,40){1.5}
%				\Dotline(15,60)(70,60){1.5}
		
		%Ponemos nombre a los ticks de los ejes
		\put(93,0){{\scriptsize$i-1$}}
		\put(128,0){{\scriptsize$i$}}
		\put(153,0){{\scriptsize$i+1$}}
		\put(20,18){{\scriptsize$j-1$}}
		\put(19,38){{\scriptsize$j$}}
		\put(0,58){{\scriptsize$j+1$}}
		
		%Ponemos nombre a los tri\~{A}!`ngulos up
		\put(90,45){{\scriptsize$T_{i-1,j}^U$}}
		\put(123,45){{\scriptsize$T_{i,j}^U$}}
		\put(110,65){{\scriptsize$T_{i,j+1}^U$}}
		
		%Ponemos nombre a los tri\~{A}!`ngulos down
		\put(84,34){{\scriptsize$T_{i-1,j-1}^D$}}
		\put(117,34){{\scriptsize$T_{i,j-1}^D$}}
		\put(77,53){{\scriptsize$T_{i-1,j}^D$}}
		\put(109,51){{\scriptsize$T_{i,j}^D$}}
		\put(137,53){{\scriptsize$T_{i+1,j}^D$}}
		\put(97,74){{\scriptsize$T_{i,j+1}^D$}}
		\put(124,74){{\scriptsize$T_{i+1,j+1}^D$}}
		
		\end{picture}
	\end{center}\vspace*{-0.3cm}
	\caption{10-point stencil for the downward triangle $T_{i,j}^D\in\mathcal{T}_{h,k}$ represented in the oblique coordinate system $\{\mathbf{e}_1', \mathbf{e}_2'\}$.}
	\label{fig:stencil:down}
\end{figure}
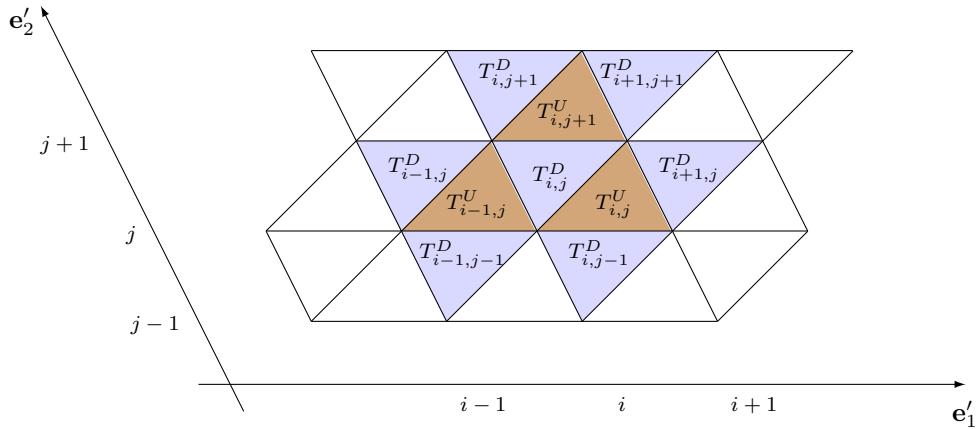

Let us consider a triangular subdomain $\Omega_k$, and let $\mathcal{T}_{h,k}\equiv\mathcal{T}_{h,k}^{\ell}$ be the conforming mesh that covers $\Omega_k$ with $4^{\ell}$ similar triangular elements. First of all, we need to set a suitable local numbering for the vertices, edges and triangles of the mesh. For a fixed refinement level $\ell$, each grid point in $\mathcal{T}_{h,k}$ is assigned a pair of indices $(i,j)$ that take the values
$j=1,2,\ldots,2^{\ell}+1$ and $i=j,j+1,\ldots,2^{\ell}+1$. Based on this notation, the vertices of $\Omega_k$ are the points $(1,1)$, $(2^{\ell}+1,1)$ and $(2^{\ell}+1,2^{\ell}+1)$. In fact, such indices can be interpreted as the coordinates of the grid points in an oblique coordinate system $\{\mathbf{e}_1',\mathbf{e}_2'\}$ that depends on the geometry of the subdomain. Further, the edges associated to the $(i,j)$-th grid node are denoted as $e_{i,j}^{m}$, where $m=1,2,3$. This latter index determines which side of the coarse triangle $\Omega_k$ is parallel to the edge under consideration. Figure \ref{fig:new:ref:sys} shows the local numbering of a mesh $\mathcal{T}_{h,k}^2$ covering $\Omega_k$ that consists of 16 triangles. The three edges associated to the grid point $(3,2)$, namely $e_{3,2}^1$, $e_{3,2}^2$ and $e_{3,2}^3$, are depicted in blue, red and brown, respectively. Note that the triangles conforming the mesh can be classified, depending on its orientation, into upward and downward triangles. As shown in Figure \ref{fig:up:down}, we will use the notation $T_{i,j}^U$ and $T_{i,j}^D$ to refer to the upward and downward triangles, respectively, that share the grid nodes $(i,j)$ and $(i+1,j+1)$. Figures \ref{fig:stencil:up} and \ref{fig:stencil:down} further show the distribution of several of these triangles on a certain region of a given mesh.

In the remaining of this section, we derive the stencil coefficients associated to an interior\footnote{In this framework, a triangle in $\mathcal{T}_{h,k}$ is said to be interior if the 10-point stencil of the pressure unknown associated to such a triangle is contained in the subdomain $\Omega_k$.} upward triangle $T_{i,j}^U$. We will assume, in the sequel, that the unit vectors normal to the edges of $\mathcal{T}_{h,k}$ are defined to point outwards from $T_{i,j}^U$ and inwards to $T_{i,j}^D$.

To begin with, let us consider $w=w_{i,j}^U$ in the equation (\ref{emfe:a}), where $w_{i,j}^U$ is the basis function of $W_{h,k}$ associated to $T_{i,j}^U$ (i.e., the characteristic function of such a triangle). If we apply the divergence theorem on $T_{i,j}^U$ and take into account the orientation of normal unit vectors on $T_{i,j}^U$, we obtain
\begin{equation}
	\label{eq:stencil:1}
	|T_{i,j}^U|\,\frac{dP_{i,j}^U}{dt}+l_2\,U_{i,j}^2+l_3\,U_{i,j}^3+l_1\,U_{i+1,j}^1=\int_{T_{i,j}^U}f\,d\mathbf{x},
\end{equation}
where $P_{i,j}^U$ denotes the pressure unknown associated to the triangle $T_{i,j}^U$. In turn, $U_{i,j}^2$, $U_{i,j}^3$ and $U_{i+1,j}^1$ stand for the components of vector $U_h$ corresponding to the edges $e_{i,j}^2$, $e_{i,j}^3$ and $e_{i+1,j}^1$, respectively. As shown in Figure \ref{fig:up:down} (left), these are the sides of $T_{i,j}^U$. Their respective lengths are denoted by $l_1$, $l_2$ and $l_3$. Note that, due to the hierarchical structure of the grid covering $\Omega_k$, the edges $e_{i,j}^m$ all have the same length $l_m$, for $m=1,2,3$, regardless of the grid point coordinates $(i,j)$.

Next, we take $\mathbf{v}=\mathbf{v}_{i,j}^2$ in the equation (\ref{semidiscrete:c}), where $\mathbf{v}_{i,j}^2$ is the basis function of $V_{h,k}$ associated to the edge $e_{i,j}^2$. Taking into account that $\mathbf{v}_{i,j}^2$ is different from zero only on the triangles $T_{i,j}^U$ and $T_{i,j}^D$ (cf. Figure \ref{fig:up:down} (left)), and using the quadrature rule introduced in (\ref{quadrature}), we get
\begin{equation}
	\label{eq:stencil:2}
	\begin{array}{rl}
		U_{i,j}^2=&\displaystyle\frac{\sqrt{3}}{l_2^2}\left(
		\tilde{U}_{i,j}^2\displaystyle\int_{T_{i,j}^U}\left(\mathbf{v}_{i,j}^{2}\right)^TGKG\,\mathbf{v}_{i,j}^{2}\,d\mathbf{x}+
		\tilde{U}_{i,j}^3\displaystyle\int_{T_{i,j}^U}\left(\mathbf{v}_{i,j}^{2}\right)^TGKG\,\mathbf{v}_{i,j}^{3}\,d\mathbf{x}\,+\right.\\[3ex]&
		\tilde{U}_{i+1,j}^1\displaystyle\int_{T_{i,j}^U}\left(\mathbf{v}_{i,j}^{2}\right)^TGKG\,\mathbf{v}_{i+1,j}^{1}\,d\mathbf{x}+
		\tilde{U}_{i,j}^2\displaystyle\int_{T_{i,j}^D}\left(\mathbf{v}_{i,j}^{2}\right)^TGKG\,\mathbf{v}_{i,j}^{2}\,d\mathbf{x}\,+\\[3ex]&
		\left.\tilde{U}_{i,j+1}^3\displaystyle\int_{T_{i,j}^D}\left(\mathbf{v}_{i,j}^{2}\right)^TGKG\,\mathbf{v}_{i,j+1}^{3}\,d\mathbf{x}+
		\tilde{U}_{i,j}^1\displaystyle\int_{T_{i,j}^D}\left(\mathbf{v}_{i,j}^{2}\right)^TGKG\,\mathbf{v}_{i,j}^{1}\,d\mathbf{x}\right),
	\end{array}
\end{equation}
where $\tilde{U}_{i,j}^m$ denotes the component of vector $\tilde{U}_h$ associated to the edge $e_{i,j}^m$. Note that the five components of vector $\tilde{U}_h$ contained in the previous formula refer to the sides of the triangles $T_{i,j}^U$ and $T_{i,j}^D$, as shown in Figure \ref{fig:up:down} (left).

In order to compute the integrals on the right-hand side of (\ref{eq:stencil:2}), we first map each of them onto the reference element $\hat{T}$. Specifically, we introduce an affine mapping $F_T$, defined from the reference triangle $\hat{T}$, with vertices $\hat{\mathbf{r}}_1=(-1,0)^T$, $\hat{\mathbf{r}}_2=(1,0)^T$ and $\hat{\mathbf{r}}_3=(0,\sqrt{3})^T$, onto any triangle $T\in\mathcal{T}_{h,k}$, with vertices $\mathbf{r}_1=(x_1,y_1)^T$, $\mathbf{r}_2=(x_2,y_2)^T$ and $\mathbf{r}_3=(x_3,y_3)^T$, i.e.,
$$
\mathbf{x}=F_T(\hat{\mathbf{x}})=B_T\,\hat{\mathbf{x}}+\mathbf{b}_T,
$$
where
$$
B_T=\left(\begin{array}{cc}\dfrac{x_2-x_1}{2}&\dfrac{2x_3-x_2-x_1}{2\sqrt{3}}\\[2ex]
\dfrac{y_2-y_1}{2}&\dfrac{2y_3-y_2-y_1}{2\sqrt{3}}\end{array}\right),\qquad \mathbf{b}_T=\left(\begin{array}{c}\dfrac{x_1+x_2}{2}\\[2ex]\dfrac{y_1+y_2}{2}\end{array}\right).
$$
Then, we consider the first integral on the right-hand side of (\ref{eq:stencil:2}) and denote $T\equiv T_{i,j}^U$. Using the fact that $(\mathbf{v}_{i,j}^{2}\circ F_T)(\hat{\mathbf{x}})=J_T^{-1}B_T\hat{\mathbf{v}}_{i,j}^{2}(\hat{\mathbf{x}})$ and $G|_T=J_TB_T^{-T}B_T^{-1}$, we obtain
$$\displaystyle\int_{T}\left(\mathbf{v}_{i,j}^{2}\right)^TGKG\,\mathbf{v}_{i,j}^{2}\,d\mathbf{x}=
J_T\displaystyle\int_{\hat{T}}\left(\hat{\mathbf{v}}_{i,j}^{2}\right)^TB_T^{-1}KB_T^{-T}\hat{\mathbf{v}}_{i,j}^{2}\,d\hat{\mathbf{x}}.
$$
This same procedure can be applied to the rest of integrals in (\ref{eq:stencil:2}). In the sequel, we show that the integral on the right-hand side of the preceding expression is independent of the particular (upward or downward) triangle under consideration.

Let us consider the local ordering of vertices displayed on Figure \ref{fig:up:down} (right) for either an upward or a downward triangle. As the spatial mesh $\mathcal{T}_{h,k}$ is composed of similar triangles, it is easy to see that all mappings from $\hat{T}$ to an upward (or downward) triangle share the same Jacobian matrix $B_{T^{U}}$ (or $B_{T^D}$). Since $B_{T^U}=-B_{T^D}$, if tensor $K$ is assumed to be constant on $\Omega_k$, then $B_{T^U}^{-1}KB_{T^U}^{-T}=B_{T^D}^{-1}KB_{T^D}^{-T}$ for all the triangles in $\mathcal{T}_{h,k}$. We denote this product by $B_{k}^{-1}KB_{k}^{-T}$, and further define $J_k\equiv J_{T^U}=J_{T^D}$. In addition, it can be shown that $\hat{\mathbf{v}}_{i,j}^m=\frac{l_m}{2}\,\hat{\mathbf{v}}_m$, for any grid point $(i,j)$ and $m=1,2,3$, where
$$\hat{\mathbf{v}}_1=\frac{1}{\sqrt{3}}\left(\begin{array}{c}1+\hat{x}\\[0ex]\hat{y}\end{array}\right),\quad \hat{\mathbf{v}}_2=\frac{1}{\sqrt{3}}\left(\begin{array}{c}-1+\hat{x}\\[0ex]\hat{y}\end{array}\right),\quad
\hat{\mathbf{v}}_3=\frac{1}{\sqrt{3}}\left(\begin{array}{c}\hat{x}\\[0ex]-\sqrt{3}+\hat{y}\end{array}\right)$$	
are the standard basis functions of $\hat{V}_h(\hat{T})$. Therefore, introducing the notation\footnote{Note that, by construction, $a_{\ell,m}=a_{m,\ell}$.}
\begin{equation}\label{a:coefficient}
a_{\ell,m}=\frac{J_k}{4}\int_{\hat{T}}\hat{\mathbf{v}}_{\ell}^TB_k^{-1}KB_k^{-T} \hat{\mathbf{v}}_m\,d\hat{\mathbf{x}},
\end{equation}
for $\ell,m\in\{1,2,3\}$, it is easy to derive, from (\ref{eq:stencil:2}), the following expression for $U_{i,j}^2$
\begin{equation}\label{eq:stencil:2bis}
	U_{i,j}^2=\frac{\sqrt{3}}{l_2}\left(2l_2\,a_{2,2}\,\tilde{U}_{i,j}^2+l_3\,a_{2,3}\left(\tilde{U}_{i,j}^3+\tilde{U}_{i+1,j}^3\right)+l_1\,a_{2,1}\left(\tilde{U}_{i,j}^1+\tilde{U}_{i+1,j}^1\right)\right).
\end{equation}

Accordingly, if we consider  $\mathbf{v}=\mathbf{v}_{i,j}^3$ in the equation (\ref{semidiscrete:c}), we will be able to express $U_{i,j}^3$ in terms of $\tilde{U}_{i,j}^3$, $\tilde{U}_{i,j-1}^1$, $\tilde{U}_{i+1,j}^1$, $\tilde{U}_{i,j}^2$ and $\tilde{U}_{i,j-1}^2$. Finally, if we take  $\mathbf{v}=\mathbf{v}_{i+1,j}^1$ in (\ref{semidiscrete:c}), $U_{i+1,j}^1$ will be given in terms of $\tilde{U}_{i+1,j}^1$, $\tilde{U}_{i,j}^2$, $\tilde{U}_{i+1,j}^2$, $\tilde{U}_{i,j}^3$ and $\tilde{U}_{i+1,j+1}^3$.

To conclude, we use the equation (\ref{semidiscrete:b}) to express the components of vector $\tilde{U}_h$ as differences of pressure unknowns. In particular, let us consider $\mathbf{v}=\mathbf{u}_{i,j}^2$ in (\ref{semidiscrete:b}). If we make use of the quadrature rule (\ref{quadrature}) and recall the orientation of normal unit vectors, we obtain
\begin{equation}
	\label{eq:stencil:3}
	\tilde{U}_{i,j}^2=\frac{\sqrt{3}}{l_2}\left(P_{i,j}^U-P_{i,j}^D\right).
\end{equation}
In this expression, the component of $\tilde{U}_h$ corresponding to the edge $e_{i,j}^2$ is written in terms of the pressure unknowns associated to the triangles that share this edge. As shown in Figure \ref{fig:up:down} (left), such triangles are $T_{i,j}^U$ and $T_{i,j}^D$. If we repeat this procedure for the remaining 8 components of $\tilde{U}_h$ involved in the expressions above, we finally obtain 10 pressure unknowns: 7 of them are associated to upward triangles and the other 3 correspond to downward triangles (depicted in brown and blue, respectively, in Figure \ref{fig:stencil:up}).

As a result, if we insert into (\ref{eq:stencil:1}) the expressions for $U_{i,j}^m$ in terms of $\tilde{U}_{i,j}^m$ derived above, for $m=1,2,3$ (see (\ref{eq:stencil:2bis}) for the case $U_{i,j}^2$), and further replace such $\tilde{U}_{i,j}^m$ by their corresponding formulas in terms of pressure differences (see (\ref{eq:stencil:3}) for the case $\tilde{U}_{i,j}^2$), we obtain an equation for the pressure variable at $T_{i,j}^U$ involving 10 unknowns, as shown in Figure \ref{fig:stencil:up}. The term $MP_h$ arising in such an equation can be explicitly written using the standard notation for the elements of a stencil
$$
S=\left(\begin{array}{lll}S_{-1,1} & S_{0,1} & S_{1,1}\\
S_{-1,0} & S_{0,0} & S_{1,0}\\
S_{-1,-1} & S_{0,-1} & S_{1,-1}\end{array}\right).
$$
More precisely,
$$
(MP_h)_{i,j}^U=\displaystyle\sum_{\ell=-1}^{1}\sum_{m=-1}^{1}S^{UU}_{\ell,m}P_{i+\ell,j+m}^U+\sum_{\ell=-1}^{1}\sum_{m=-1}^{1}S^{UD}_{\ell,m}P_{i+\ell,j+m}^D,
$$
where
$$
S^{UU}=3\left(\begin{array}{ccc} 0 & a_{2,3} & a_{1,3}\\[1ex]
a_{1,2} & 2\displaystyle\sum_{\ell=1}^3\sum_{m=\ell}^3a_{\ell,m} & a_{1,2}\\[1ex]
a_{1,3} & a_{2,3} & 0\end{array}\right),\
S^{UD}=-6\left(\begin{array}{ccc}0 & 0 & 0\\[1ex]
0 & \displaystyle\sum_{\ell=1}^3 a_{\ell,2} & \displaystyle\sum_{\ell=1}^3 a_{1,\ell}\\[1ex]
0 & \displaystyle\sum_{\ell=1}^3 a_{\ell,3} & 0 \end{array}\right).
$$
The entries $a_{\ell,m}$ in both matrices are given by (\ref{a:coefficient}). Remarkably, the stencil coefficients are the same for all interior pressures associated to upward triangles. Furthermore, these coefficients depend on tensor $K$ and the geometry of the given subdomain $\Omega_k$, but not on the number of mesh refinement processes.

The procedure to derive the stencil coefficients associated to an interior downward triangle is completely analogous. In such a case, we obtain
$$
(MP_h)_{i,j}^D=\displaystyle\sum_{\ell=-1}^{1}\sum_{m=-1}^{1}S^{DU}_{\ell,m}P_{i+\ell,j+m}^U+\sum_{\ell=-1}^{1}\sum_{m=-1}^{1}S^{DD}_{\ell,m}P_{i+\ell,j+m}^D,
$$
where $S^{DD}=S^{UU}$, and
$$
S^{DU}=-6\left(\begin{array}{ccc}0 & \displaystyle\sum_{\ell=1}^3 a_{\ell,3} & 0\\[1ex]
\displaystyle\sum_{\ell=1}^3 a_{1,\ell} & \displaystyle\sum_{\ell=1}^3 a_{\ell,2} & 0\\[1ex]
0 & 0 & 0 \end{array}\right).
$$
Figure \ref{fig:stencil:down} shows the resulting stencil involving 10 pressure unknowns: 7 of them associated to downward triangles and the other 3 corresponding to upward triangles, coloured in blue and brown, respectively.

Note that, if we consider non-interior triangles, the coefficients for the pressure stencil can be derived following the same strategy; the only difference is that, in this case, the matrices $S^{UU}$, $S^{UD}$, $S^{DU}$ and $S^{DD}$ may contain a greater number of null elements.

\vspace*{0.5cm}

\textbf{Acknowledgments.} This work was partially supported by MINECO grant MTM2014-52859.


\begin{thebibliography}{99}
\footnotesize

\bibitem{arb:daw:kee:94} T. Arbogast, C.N. Dawson, P.T. Keenan, Mixed finite elements as finite differences for elliptic equations on triangular elements, Tech. Report CRPC-TR94392, Rice University, 1994.

\bibitem{arb:daw:kee:whe:yot:98} T. Arbogast, C.N. Dawson, P.T. Keenan, M.F. Wheeler, I. Yotov, Enhanced cell-centered finite differences for elliptic equations on general geometry, SIAM J. Sci. Comput. 19 (1998) 404--425.

%\bibitem{arb:whe:yot:97} T. Arbogast, M.F. Wheeler, I. Yotov, Mixed finite elements for elliptic problems with tensor coefficients as cell-centered finite differences, SIAM J. Numer. Anal. 34 (1997) 828--852.
    
\bibitem{arn:bre:85} D.N. Arnold, F. Brezzi, Mixed and nonconforming finite element methods : implementation, postprocessing and error estimates, RAIRO Mod\'el. Math. Anal. Num\'er. 19 (1985) 7--32.
    
\bibitem{arr:por:14} A. Arrar\'as, L. Portero, Expanded mixed finite element domain decomposition methods on triangular grids, Int. J. Numer. Anal. Model. 11 (2014) 255--270.

\bibitem{ber:hul:04} B.K. Bergen, F. H$\ddot{\mathrm{u}}$lsemann, Hierarchical hybrid grids: data structures and core algorithms for multigrid, Numer. Linear Algebra Appl. 11 (2004) 279--291.

\bibitem{ber:wel:hul:rud:07} B. Bergen, G. Wellein, F. H$\ddot{\mathrm{u}}$lsemann, U. R$\ddot{\mathrm{u}}$de, Hierarchical hybrid grids: achieving TERAFLOP performance on large scale finite element simulations, Int. J. Parallel Emergent Distrib. Syst. 22 (2007) 311--329.

\bibitem{bof:bre:for:13} D. Boffi, F. Brezzi, M. Fortin, Mixed Finite Element Methods and Applications, Springer Ser. Comput. Math. 44, Springer-Verlag, Berlin, 2013.

\bibitem{bre:for:mar:06} F. Brezzi, M. Fortin, L.D. Marini, Error analysis of piecewise constant pressure approximations of Darcy's law, Comput. Methods Appl. Mech. Engrg. 195 (2006) 1547--1559.

\bibitem{che:98} Z. Chen, Expanded mixed finite element methods for linear second-order elliptic problems I, RAIRO Mod\'el. Math. Anal. Num\'er. 32 (1998) 479--499.

\bibitem{dup:kee:98} T.F. Dupont, P.T. Keenan, Superconvergence and postprocessing of fluxes from lowest-order mixed methods on triangles and tetrahedra, SIAM J. Sci. Comput. 19 (1998) 1322--1332.

\bibitem{ewi:sae:she:98} R.E. Ewing, O. S{\ae}vareid, J. Shen, Discretization schemes on triangular grids, Comput. Methods Appl. Mech. Engrg. 152 (1998) 219--238.

\bibitem{gan:lis:sha:zen:02} V. Ganzha, R. Liska, M. Shashkov, C. Zenger, Support operator method for Laplace equation on unstructured triangular grid, Sel\c{c}uk J. Appl. Math. 3 (2002) 21--48.

\bibitem{gat:heu:01} G.N. Gatica, N. Heuer, An expanded mixed finite element approach via a dual--dual formulation and the minimum residual method, J. Comput. Appl. Math. 132 (2001) 371--385.

\bibitem{gme:rud:14} B. Gmeiner, U. R$\ddot{\mathrm{u}}$de, Peta-scale hierarchical hybrid multigrid using hybrid parallelization, in: Large-scale scientific computing, I. Lirkov, S. Margenov, J. Wa$\acute{\mathrm{s}}$niewski (eds.), Lecture Notes in Comput. Sci. 8353, Springer-Verlag, Berlin, 2014, pp. 439--447.

\bibitem{gme:rud:ste:wal:woh:15} B. Gmeiner, U. R$\ddot{\mathrm{u}}$de, H. Stengel, C. Waluga, B. Wohlmuth, Performance and scalability of hierarchical hybrid multigrid solvers for Stokes systems, SIAM J. Sci. Comput. 37 (2015) C143--C168.

\bibitem{hai:lub:roc:89} E. Hairer, Ch. Lubich, M. Roche, The Numerical Solution of Differential-Algebraic Systems by Runge--Kutta Methods, Lecture Notes in Math. 1409, Springer-Verlag, Berlin, 1989.

\bibitem{hai:wan:96} E. Hairer, G. Wanner, Solving Ordinary Differential Equations II. Stiff and Differential-Algebraic Problems, Second edition, Springer Ser. Comput. Math., Vol. 14, Springer-Verlag, Berlin, 1996.

\bibitem{mac:car:88} R.J. MacKinnon, G.F. Carey, Analysis of material interface discontinuities and superconvergent fluxes in finite difference theory, J. Comput. Phys. 75 (1988) 151--167.

\bibitem{mat:08} T.P.A. Mathew, Domain Decomposition Methods for the Numerical Solution of Partial Differential Equations, Lect. Notes Comput. Sci. Eng. 61, Springer-Verlag, Berlin, 2008.

\bibitem{mic:sac:sal:01} S. Micheletti, R. Sacco, F. Saleri, On some mixed finite element methods with numerical integration, SIAM J. Sci. Comput. 23 (2001) 245--270.

\bibitem{rav:tho:77} R.A. Raviart, J.M. Thomas, A mixed finite element method for 2nd order elliptic problems, in: Mathematical Aspects of Finite Element Methods, I. Galligani, E. Magenes (eds.), Lecture Notes in Math. 606, Springer, Berlin, 1977, pp. 292--315.

\bibitem{rod:gas:lis:12} C. Rodrigo, F.J. Gaspar, F.J. Lisbona, Multigrid methods on semi-structured grids, Arch. Comput. Methods Eng. 19 (2012) 499--538.

\bibitem{tho:77} J.M. Thomas, Sur l'analyse num\'erique des m\'ethods d'\'el\'ements finis hybrides et mixtes, Ph.D. thesis, Universit\'e Pierre et Marie Curie, 1977.

\bibitem{tho:97} V. Thom\'ee, Galerkin Finite Element Methods for Parabolic Problems, Springer Ser. Comput. Math. 25, Springer-Verlag, Berlin, 1997.

\bibitem{whe:73} M.F. Wheeler, A priori ${L}_2$ error estimates for Galerkin approximations to parabolic partial differential equations, SIAM J. Numer. Anal. 10 (1973) 723--759.

\end{thebibliography}
\end{document}